\documentclass{amsart}[12pt]

\usepackage{enumerate}
\usepackage[hidelinks]{hyperref}
\usepackage{verbatim}
\usepackage{yfonts} %
\usepackage{amssymb} %
\usepackage{amsthm}
\usepackage{array}
\usepackage{booktabs}%
\usepackage{hhline}%
\usepackage{xy} %
\usepackage{epsfig}%
\usepackage{color}%
\usepackage{upgreek}
\usepackage[english]{babel}
\usepackage{epigraph}%
\usepackage{fancybox}%
\usepackage{shadow}
\usepackage{afterpage}
\usepackage{mathrsfs}
\usepackage{enumitem}
\usepackage{tabularx}
\usepackage{subcaption}
\usepackage{graphicx}
\usepackage{type1cm}
\usepackage{eso-pic}
\usepackage{color}
\usepackage{upgreek}
\usepackage{bigints}
\usepackage{thmtools}
\usepackage{thm-restate}

\newtheorem{theorem}{Theorem}
\newtheorem{proposition}[theorem]{Proposition}

\newtheorem{lemma}[theorem]{Lemma}

\theoremstyle{definition}
\newtheorem{definition}[theorem]{Definition}
\newtheorem{remark}[theorem]{Remark}
\newtheorem{example}[theorem]{Example}

\theoremstyle{plain}

\newcommand{\vt}{\vspace{.1cm}}

\newcommand{\R}{\mathbb{R} }
\newcommand{\q}{\mathbb{Q} }
\newcommand{\qq}{\q_{\epsilon_1}^{n_1}\times\q_{\epsilon_2}^{n_2}}

\newcommand{\h}{\mathbb{H}}

\newcommand{\s}{\mathbb{S}}
\newcommand{\sone}{s_{\tau_1}}
\newcommand{\stwo}{s_{\tau_2}}
\newcommand{\cone}{c_{\tau_1}}
\newcommand{\ctwo}{c_{\tau_2}}

\renewcommand{\rho}{\varrho}
\renewcommand{\theta}{\varTheta}
\renewcommand{\Theta}{\varTheta}
\renewcommand{\Sigma}{\varSigma}
\renewcommand{\Omega}{\varOmega}
\renewcommand{\Lambda}{\varLambda}
\renewcommand{\emptyset}{\varnothing}
\renewcommand{\tau}{\uptau}
\usepackage{amsmath}

\newcommand{\overbar}[1]{\mkern 1.5mu\overline{\mkern-1.5mu#1\mkern-1.5mu}\mkern 1.5mu}

\newcommand{\la}{\langle}
\newcommand{\ra}{\rangle}

\makeatletter
\newcommand{\tpitchfork}{%
  \vbox{
    \baselineskip\z@skip
    \lineskip-.52ex
    \lineskiplimit\maxdimen
    \m@th
    \ialign{##\crcr\hidewidth\smash{$-$}\hidewidth\crcr$\pitchfork$\crcr}
  }%
}
\makeatother

\DeclareMathOperator{\spec}{spec}
\DeclareMathOperator{\tr}{tr}

\begin{document}

\title[]
{Isoparametric Hypersurfaces in Products of
Simply Connected Space Forms}
\author{Ronaldo F. de Lima \and Giuseppe Pipoli}

\address[A1]{Departamento de Matem\'atica - Universidade Federal do Rio Grande do Norte}
\email{ronaldo.freire@ufrn.br}
\address[A2]{Department of Information Engineering, Computer Science and Mathematics, Università degli Studi dell’Aquila.}
\email{giuseppe.pipoli@univaq.it}

\maketitle

\begin{abstract}
For $i\in\{1,2\}, $ let $\q_{\epsilon_i}^{n_i}$ denote the simply connected
space form of dimension $n_i\ge 2$ and
constant sectional curvature $\epsilon_i\in\R$,
$(\epsilon_1,\epsilon_2)\ne (0,0)$.
We prove that any connected isoparametric hypersurface  of
\,$\qq$ has constant angle function.
We use this property to
classify both the isoparametric and homogeneous
hypersurfaces of $\qq$ when $\epsilon_1\epsilon_2\le 0$. Under additional hypotheses,
we obtain a similar classification result in the case $\epsilon_1\epsilon_2> 0$.

\vspace{.3cm}

\vspace{.1cm}

\end{abstract}

\section{Introduction}
The isoparametric hypersurfaces of a Riemannian manifold $M^n$
are those whose nearby parallel hypersurfaces, together with the hypersurface itself,
have constant mean curvature. They arise as level sets of smooth functions
$F$ on $M^n$, also called isoparametric, whose gradient $\nabla F$ and Laplacian
$\Delta F$ satisfy
$$
\|\nabla F\|^2=\Phi\circ F, \qquad \Delta F=\Psi\circ F,
$$
where $\Phi$ and $\Psi$ are smooth real-valued functions.

The classification of isoparametric hypersurfaces is a major topic in
submanifold theory. The first results in this direction were obtained
in the late 1930s by Levi-Civita~\cite{L-C}, who classified the
isoparametric surfaces of $\R^3$; Segre~\cite{segre}, who extended
Levi-Civita's result to the Euclidean spaces $\R^n$; and
Cartan~\cite{cartan1,cartan3}, who classified the isoparametric
hypersurfaces of the hyperbolic spaces $\h^n$.

In~\cite{cartan2,cartan4}, Cartan obtained many results on isoparametric
hypersurfaces of the spheres $\s^n$, including the classification of those
having at most three distinct principal curvatures.
In fact, the classification problem for isoparametric hypersurfaces of
the spheres $\s^n$ is rather involved, and a complete solution has only
recently been announced; see~\cite{cecil,chi,siffert1,siffert2}.

The homogeneous hypersurfaces of a Riemannian manifold $M^n$ --- that is,
the codimension-one orbits of actions of Lie subgroups of ${\rm Iso}(M)$ ---
are isoparametric and have constant principal curvatures.
The converse, however, does not hold in general. In fact, one of the major
problems in this theory is to determine when isoparametric hypersurfaces
are necessarily homogeneous or have constant principal curvatures;
cf.~\cite{DV}.

In his aforementioned works, Cartan showed that, when the ambient manifold
$M^n$ is a simply connected space form, being isoparametric is equivalent
to having constant principal curvatures. Moreover, these conditions are
equivalent to homogeneity when $M^n$ is either $\R^n$ or $\h^n$.
On the other hand, spheres $\s^n$ admit nonhomogeneous isoparametric
hypersurfaces, such as the FKM examples; see~\cite{cecil,FKM}.
The classification of homogeneous hypersurfaces of $\s^n$ is due to
Hsiang and Lawson~\cite{HL}.

In this paper, we investigate isoparametric hypersurfaces of the product spaces
$\qq$, where $\q_{\epsilon_i}^{n_i}$ denotes the simply connected
space form of dimension $n_i\ge 2$ and constant sectional curvature
$\epsilon_i\in\R$. Throughout the paper, we assume that
these spaces are non-Euclidean, that is,
$$
(\epsilon_1,\epsilon_2)\ne(0,0).
$$

On every hypersurface of $\qq$, there is a natural \emph{angle function} $\theta$
whose properties are similar to those of the angle function of hypersurfaces in
cylinders $\q_\epsilon^n\times\R$; see Sect.~\ref{sec-constantangle}.
For instance, hypersurfaces with constant angle $\theta=\pm1$ are of one of the
following two types:
$\Sigma_1^{n_1-1}\times\q_{\epsilon_2}^{n_2}$ or
$\q_{\epsilon_1}^{n_1}\times\Sigma_2^{n_2-1}$,
where $\Sigma_i^{n_i-1}$ is a hypersurface of
$\q_{\epsilon_i}^{n_i}$; see Prop.~\ref{prop-Ni=0}.
In~\cite{dLP}, we proved that every isoparametric hypersurface of
$\q_\epsilon^n\times\R$ has constant angle. Our first result shows that
the same holds for the products $\qq$.

\begin{restatable}{theorem}{isoimpliesthetact}
\label{thm-isoparametric-->thetaconstant}
Every connected isoparametric hypersurface of the product \,$\qq$
has constant angle function.
\end{restatable}

The proof of Theorem~\ref{thm-isoparametric-->thetaconstant}
follows the same lines as that of~\cite[Prop.~7]{dLP}:
by applying Jacobi field theory, we reduce it
to the solution of an intricate algebraic problem; see
Sect.~\ref{sec-proofTh1}. However, the argument we give here
is different, and can be applied to provide a simpler proof
of~\cite[Prop.~7]{dLP}.
We add that the condition $(\epsilon_1,\epsilon_2)\ne(0,0)$ is necessary
in Theorem~\ref{thm-isoparametric-->thetaconstant}, since round spheres in Euclidean space are isoparametric
hypersurfaces with nonconstant angle function.

We apply Theorem~\ref{thm-isoparametric-->thetaconstant}
to classify both the isoparametric and the homogeneous hypersurfaces
of $\qq$ when $\epsilon_1\epsilon_2\le 0$. The result reads as follows.

\begin{restatable}{theorem}{classificationa}
\label{thm-classification01}
Assume $\epsilon_1\epsilon_2\le 0$ and let $\Sigma\subset\qq$ be a connected oriented hypersurface.
Then the following conditions are equivalent:
\begin{enumerate}[parsep=1ex]
 \item[\textup{(i)}] $\Sigma$ is isoparametric.

 \item[\textup{(ii)}] $\Sigma$ is an open set of one of the following
 complete hypersurfaces:
 \begin{itemize}[parsep=1ex]
  \item[\textup{(a)}]
  $\Sigma_1^{n_1-1}\times\q_{\epsilon_2}^{n_2}$, where
  $\Sigma_1^{n_1-1}\subset\q_{\epsilon_1}^{n_1}$ is complete and
  isoparametric;
   \item[\textup{(b)}]
  $\q_{\epsilon_1}^{n_1}\times\Sigma_2^{n_2-1}$, where
  $\Sigma_2^{n_2-1}\subset\q_{\epsilon_2}^{n_2}$ is complete and
  isoparametric;
  \item[\textup{(c)}]
  a flat-horospherical hypersurface, when $\epsilon_1, \epsilon_2\le 0.$
 \end{itemize}
\end{enumerate}
In particular, every isoparametric hypersurface of \,$\qq$ has constant
principal curvatures.
Moreover, if $\epsilon_1,\epsilon_2\leq0$, conditions \textup{(i)} and
\textup{(ii)} are also equivalent to
\begin{enumerate}
 \item[\textup{(iii)}]
 $\Sigma$ is an open set of a homogeneous hypersurface.
\end{enumerate}
If one factor is spherical, i.e., if $\epsilon_i>0$ for some
$i\in\{1,2\}$, then condition \textup{(iii)} is equivalent to
\begin{enumerate}
 \item[\textup{(iv)}]
 $\Sigma$ is an open set of a hypersurface of type \textup{(a)} or
 \textup{(b)}, and the corresponding hypersurface
 $\Sigma_i^{n_i-1}\subset\q_{\epsilon_i}^{n_i}$ is homogeneous.
\end{enumerate}
\end{restatable}

\begin{figure}[hbt]
\includegraphics[width=10cm]{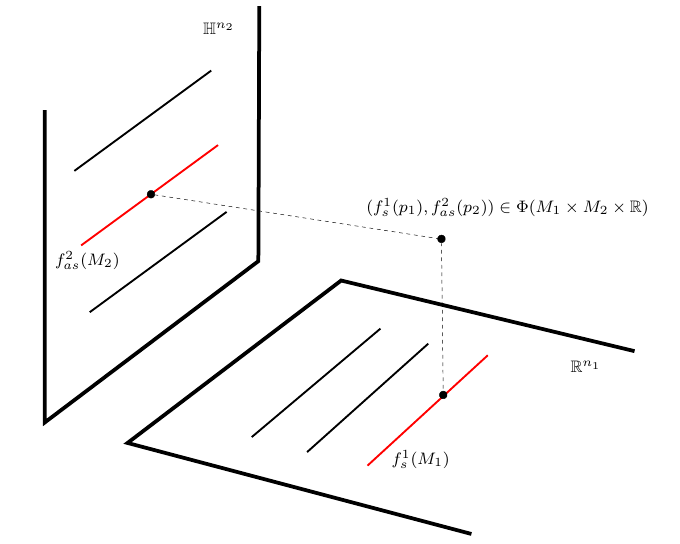}
\caption{\small Coordinates on an immersed  flat-horospherical hypersurface
$\Phi(M_1\times M_2\times\R)\subset\R^{n_1}\times\h^{n_2}_{\epsilon_2}$,
$\Phi(p_1,p_2, s):=(f_s^1(p_1), f_{as}^2(p_2))$, where $a>0$ and $\{f_s^1\}$ (resp.  $\{f_{as}^2\}$)
is a parallel family of hyperplanes of $\R^{n_1}$ (resp.~horospheres of $\h^{n_2}_{\epsilon_2}$);
see Ex.~\ref{examp-horospherical}, Sect.~\ref{sec-constantangle}.}
\label{fig-flathorospherical}
\end{figure}

The hypersurface in item \textup{(ii)(c)} of
Theorem~\ref{thm-classification01}
is built from two one-parameter families $f_s^i$
of parallel hypersurfaces in $\q_{\epsilon_i}^{n_i}$, $\epsilon_i\le 0$,
which are related by a linear function $\phi(s)=as$, $s\in\R$.
One family consists of parallel hyperplanes in the Euclidean space
$\R^{n_i}$, whereas the other consists of parallel horospheres in the
hyperbolic space $\h_{\epsilon_j}^{n_j}$ of constant curvature
$\epsilon_j<0$, $\{i,j\}=\{1,2\}$.
The resulting hypersurface is called \emph{flat-horospherical}
(Fig.~\ref{fig-flathorospherical}).
In an analogous way, we construct isoparametric hypersurfaces in
$\h_{\epsilon_1}^{n_1}\times\h_{\epsilon_2}^{n_2}$,
called \emph{bi-horospherical}, from two families of parallel
horospheres, one in each factor; see Ex.~\ref{examp-horospherical},  Sect.~\ref{sec-constantangle}.

The classification of isoparametric hypersurfaces of $\qq$ when
$\epsilon_1\epsilon_2>0$ is a more involved problem. In this case,
we obtain a partial solution under two additional hypotheses.

First, we assume that the hypersurface has a \emph{distinguished point}
$p$. This means that either $\theta^2(p)=1$ or the
principal directions of $\Sigma$ at  $p$ \emph{split properly},
i.e., $n_i-1$ of them are tangent to $\q_{\epsilon_i}^{n_i}$,
$i=1,2$; see Defs.~\ref{def-split} and~\ref{def:distinguished}.

Second, assuming that $\theta^2<1$, we define on the hypersurface the
functions
$$
\tau_i=-\epsilon_i(1\pm\theta)/2.
$$
When the pairs $(\sqrt{\tau_1(p)},\sqrt{\tau_2(p)})\in\mathbb C^2$
and $(n_1,n_2)\in\mathbb N^2$ satisfy a certain relation, we say that
the former is $(n_1,n_2)$-dependent. Otherwise, it is said to be
$(n_1,n_2)$-independent; see Def.~\ref{n1n2dep}. So, our second hypothesis
is that $(\sqrt{\tau_1(p)},\sqrt{\tau_2(p)})$ is
$(n_1,n_2)$-independent.

\begin{restatable}{theorem}{classificationb}
\label{thm-classification02}
Suppose that  $\epsilon_1\epsilon_2> 0$, and let $\Sigma\subset\qq$ be a connected oriented hypersurface
with a distinguished point $p$. If~$\Theta(p)\ne\pm 1$, assume in addition that
$(\sqrt{\tau_1(p)},\sqrt{\tau_2(p)})$ is $(n_1,n_2)$-independent.
Then the following conditions are equivalent:
\begin{enumerate}[parsep=1ex]
 \item[\textup{(i)}]
 $\Sigma$ is isoparametric.

 \item[\textup{(ii)}]
 $\Sigma$ is an open set of one of the following complete
 hypersurfaces:
 \begin{itemize}[parsep=1ex]
  \item[\textup{(a)}]$\Sigma_1^{n_1-1}\times\q_{\epsilon_2}^{n_2}$, where
  $\Sigma_1^{n_1-1}\subset\q_{\epsilon_1}^{n_1}$ is complete and
  isoparametric;

  \item[\textup{(b)}]
  $\q_{\epsilon_1}^{n_1}\times\Sigma_2^{n_2-1}$, where
  $\Sigma_2^{n_2-1}\subset\q_{\epsilon_2}^{n_2}$ is complete and
  isoparametric;

  \item[\textup{(c)}]
  a bi-horospherical hypersurface, if $\epsilon_1,\epsilon_2<0$.
 \end{itemize}
\end{enumerate}
In particular, every isoparametric hypersurface of \,$\qq$ satisfying
the preceding assumptions has constant principal curvatures.
Moreover, if both $\epsilon_1$ and $\epsilon_2$ are negative, conditions \textup{(i)} and
\textup{(ii)} are also equivalent to
\begin{enumerate}
 \item[\textup{(iii)}]
 $\Sigma$ is an open set of a homogeneous hypersurface.
\end{enumerate}
If $\epsilon_1,\epsilon_2>0$, then condition \textup{(iii)} is
equivalent to
\begin{enumerate}
 \item[\textup{(iv)}]
 $\Sigma$ is an open set of a hypersurface of type \textup{(a)} or
 \textup{(b)}, and the corresponding hypersurface
 $\Sigma_i^{n_i-1}\subset\s_{\epsilon_i}^{n_i}$ is homogeneous.
\end{enumerate}
\end{restatable}

We note that  there exist
one-parameter families of homogeneous hypersurfaces
in $\h^n\times\h^n$ and in $\s^n\times\s^n$
with no distinguished points; cf. Ex.~\ref{exam-tube},~\ref{exam-cui}. This shows
that the existence of a distinguished point is a necessary hypothesis in
Theorem~\ref{thm-classification02}.

The isoparametric (resp.~homogeneous) hypersurfaces of
$\q_{\epsilon_1}^2\times\q_{\epsilon_2}^2$, $\epsilon_1\ne\epsilon_2$,
have been classified in~\cite{gao01} (resp.~\cite{santos}), and
the isoparametric and  homogeneous  hypersurfaces of
$\h^2\times\h^2$ (resp.~$\s^2\times\s^2$) have been classified
in~\cite{gao} (resp.~\cite{urbano}). In all of these works, it is verified
that  isoparametric/homogeneous hypersurfaces
have constant angle function.

The paper is organized as follows. In Section~\ref{sec-preliminaries}, we set some notation
and  discuss  a few geometric aspects of $\qq$. In Section~\ref{sec-constantangle}, we establish
results on constant angle hypersurfaces of $\qq$ that will be used in the proof of our main results.
Finally, in Sections~\ref{sec-proofTh1}--\ref{sec-classification02} we prove
Theorems~\ref{thm-isoparametric-->thetaconstant}--\ref{thm-classification02}.

\medskip
\noindent
\emph{Added in proof}. After a preliminary version of this paper had been completed,
we became aware of the work of Tan, Xie, and Yan~\cite{tan02}, whose results overlap
with ours in the case where one of the factors is Euclidean.
We note, however, that their classification results, in contrast to ours,
require the hypersurfaces to be complete.
In addition, a careful analysis of their manuscript revealed several
issues in the proofs of their main theorems; see Remarks~\ref{rem:gaps01} and~\ref{rem:gaps02} for details.
For these reasons, we do not regard the results in~\cite{tan02} as established.

\section{Preliminaries}  \label{sec-preliminaries}
\subsection{Some geometric aspects of $\qq$}
Given integers $n_1, n_2\ge 2$ and $\epsilon_1,\epsilon_2\in\R$, we shall
consider the  product space $\qq$ endowed with its standard product metric,
which will  be denoted by $\la\,,\,\ra$.

The tangent bundle $\mathfrak X(\qq)$ of $\qq$ splits as
$$\mathfrak X(\qq)=\mathfrak{X}(\q_{\epsilon_1}^{n_1})\oplus\mathfrak{X}(\q_{\epsilon_2}^{n_2}),$$
where $\mathfrak{X}(\q_{\epsilon_i}^{n_i})$ denotes the tangent bundle of $\q_{\epsilon_i}^{n_i}$.
Considering this fact, for a given
$X\in\mathfrak{X}(\q_{\epsilon_1}^{n_1})\oplus\mathfrak{X}(\q_{\epsilon_2}^{n_2})$,
we shall write
\[
X=X_1+X_2, \,\,\, X_1\in\mathfrak{X}(\q_{\epsilon_1}^{n_1})\oplus\{0\},
\,\,\, X_2\in\{0\}\oplus\mathfrak{X}(\q_{\epsilon_2}^{n_2}).
\]

The  curvature tensor of $\qq$ will be denoted by $\overbar R$, so that
\[
\overbar R(X,Y)Z=R_{\epsilon_1}(X_1,Y_1)Z_1+R_{\epsilon_2}(X_2,Y_2)Z_2,
\]
where $R_{\epsilon_i}$ is the curvature tensor of $\q_{\epsilon_i}^{n_i}$, that is,
\[
\overbar R_{\epsilon_i}(X_i,Y_i)Z_i=\epsilon_i(\la X_i,Z_i\ra Y_i-\la Y_i,Z_i\ra X_i).
\]

The operator $P$  on $\mathfrak X(\q_{\epsilon_1}^{n_1})\oplus\mathfrak X(\q_{\epsilon_2}^{n_2})$ defined by
\[
P(X_1+X_2)=X_1-X_2
\]
will play a fundamental role here. It is easily checked that $P$ is parallel, that is,
\[
\overbar\nabla_XPY=P\overbar\nabla_XY\,\,\,
\forall X,Y\in\mathfrak X(\q_{\epsilon_1}^{n_1})\oplus\mathfrak X(\q_{\epsilon_2}^{n_2}),
\]
where  $\overbar\nabla$ denotes the Levi-Civita connection of $\qq$.

\subsection{Hypersurfaces of \,$\qq$}
Let $\Sigma\subset\qq$ be an oriented hypersurface with unit normal
$N=N_1+N_2$ and shape operator $A$ with respect to $N$, that is,
\[
AX=-\overbar\nabla_XN,  \,\,\,  X\in \mathfrak X(\Sigma),
\]
where $\mathfrak X(\Sigma)$ is the tangent bundle of $\Sigma$.

The principal curvatures of $\Sigma,$ that is, the eigenvalues of the shape operator
$A,$ will be denoted by $k_1\,, \dots ,k_n$. In this setting,
we define the non-normalized \emph{mean curvature} $H$ of  $\Sigma$
as the sum of its principal curvatures:
$$H=k_1+\cdots +k_n.$$

We shall denote by $\nabla$ the Levi-Civita connection of $\Sigma$. The gradient
of a smooth function $f$ on $\Sigma$ will be denoted by $\nabla f$.

\begin{definition} \label{def-angle}
With the above notation, the \emph{angle function} $\theta$ of $\Sigma$ is defined by
\[
\theta=\langle PN,N\rangle=\|N_1\|^2-\|N_2\|^2.
\]
The tangential component of $PN$ will be denoted by $T$, so that
\[
T=PN-\theta N=(1-\theta)N_1-(1+\theta)N_2.
\]
\end{definition}

The function $\theta$ satisfies
\begin{equation} \label{eq-N1N2theta}
\theta=2\|N_1\|^2-1=-2\|N_2\|^2+1,
\end{equation}
from which we conclude that:
\begin{itemize}[parsep=1ex]
\item  $\theta$ is constant  if and only if $\|N_1\|$ and $\|N_2\|$ are both  constant;
\item  either $\|N_1\|$ or $\|N_2\|$ vanishes on $\Sigma$ if and only if $\theta=\pm 1.$
\end{itemize}

The following lemma follows from straightforward computations
(see~\cite{gao,urbano}).

\begin{lemma} \label{lem-T&theta}
The following equalities hold:
\begin{itemize}[parsep=1ex]
\item $\nabla\theta=-2AT$;
\item $\nabla_XT=\theta AX-PAX+\la PAX,N\ra N$ for all $X\in\mathfrak X(\Sigma)$.
\end{itemize}
In particular, if $T$ never vanishes and $\theta$ is constant on $\Sigma$, one has:
\begin{itemize}[parsep=1ex]
\item $T$ is a principal direction  with principal curvature $0$;
\item $\overbar\nabla_TT=\nabla_TT=0$, that is, the trajectories of \,$T$ are geodesics in both
$\Sigma$ and $\qq$.
\end{itemize}
\end{lemma}

To conclude this section, we recall the Codazzi equation
for hypersurfaces, which in our context reads as
\begin{equation} \label{eq-codazzi}
(\nabla A_X)Y-(\nabla A_Y)X=(\overbar R(X,Y)N)^\top \,\,\,\forall X,Y\in\mathfrak X(\Sigma),
\end{equation}
where $(\nabla A_X)Y=\nabla A_XY-A\nabla_XY$, and
\[
Z^\top=Z-\langle Z,N\rangle N, \,\,\, Z\in\mathfrak{X}(\q_{\epsilon_1}^{n_1})\oplus\mathfrak{X}(\q_{\epsilon_2}^{n_2}).
\]

\section{Hypersurfaces of $\qq$ with constant angle} \label{sec-constantangle}

In this section, we establish some  results on
constant angle hypersurfaces of $\qq$ that will be used in the proofs
of our main theorems.

Let us first consider hypersurfaces  of
$\qq$ which are defined as immersions $\Phi$ of the following type:
\begin{equation} \label{eq-parametrization}
\begin{array}{cccc}
\Phi\colon & M_1^{n_1-1}\times M_2^{n_2-1}\times I & \rightarrow & \qq\\
           &   (p_1,p_2,s)  & \mapsto     & (\exp_{p_1}(s\eta_{p_1}), \exp_{p_2}(\phi(s)\eta_{p_2})),
\end{array}
\end{equation}
where $I\subset\R$ is an open interval, $\phi\colon I\to\phi(I)$ is a diffeomorphism,
$M_i^{n_i-1}\hookrightarrow\q_{\epsilon_i}^{n_i}$ is a hypersurface of $\q_{\epsilon_i}^{n_i}$, and
$\eta_{p_i}$  is the unit normal to $M_i^{n_i-1}$ at $p_i$, $i\in\{1,2\}$.

Defining $f_s^i:M_i^{n_i-1}\to\q_{\epsilon_i}^{n_i}$ by
\begin{equation} \label{eq-families}
f_s^1(p_1)=\gamma_{p_1}(s):=\exp_{p_1}(s\eta_{p_1}),
\quad f_s^2(p_2)=\gamma_{p_2}(\phi(s))=:\exp_{p_2}(\phi(s)\eta_{p_2}),
\end{equation}
we have that $\{f_s^i\,;\, s\in I\}$ is a family
of parallel hypersurfaces of $\q_{\epsilon_i}^{n_i}$.

\begin{definition}
We shall call $\Sigma:=\Phi(M_1^{n_1-1}\times M_2^{n_2-1}\times I)\subset\qq$
as above an $(f_s^1,f_{\phi(s)}^2)$-\emph{hypersurface} of $\qq$.
\end{definition}

Straightforward calculations give that, at $p=(p_1,p_2,s)\in M_1^{n_1-1}\times M_2^{n_2-1}\times I,$
the unit normal $N$ to the $(f_s^1,f_{\phi(s)}^2)$-hypersurface $\Sigma$ is
\begin{equation} \label{eq-normal}
N=\frac{1}{\mathcal W(s)}\left(-\phi'(s)\gamma_{p_1}'(s)+\gamma_{p_2}'(\phi(s))\right),  \,\,\, \mathcal W(s)=\sqrt{1+(\phi'(s))^2},
\end{equation}
and that its $T$-field is
$$T=PN-\theta N=-\frac{2\phi'(s)}{\mathcal W^3(s)}\left(\gamma_{p_1}'(s)+\phi'(s)\gamma_{p_2}'(\phi(s))\right).$$

Set $n=n_1+n_2$. Denoting by
$$X^1, \dots, X^{n_1-1}\in \mathfrak X(\q_{\epsilon_1}^{n_1}), \quad
Y^1, \dots, Y^{n_2-1}\in \mathfrak X(\q_{\epsilon_2}^{n_2})$$
the unit principal directions of $M_1$ and $M_2$, respectively,
it is easily checked that the unit principal directions
$Z^1,\dots, Z^{n-1}$ of $\Sigma$ are:
\begin{align} \label{eq-principaldirections}
Z^{1\phantom{+i}}       &=T/\|T\| ;\nonumber\\[1ex]
Z^{1+i}   &=X_i\in\mathfrak{X}(\q_{\epsilon_1}^{n_1})\oplus\{0\} \,\,\, \forall i\in\{1,\dots, n_1-1\}; \\[1ex]
Z^{n_1+i}&=Y_i\in\{0\}\oplus\mathfrak{X}(\q_{\epsilon_2}^{n_2}) \,\,\, \forall i\in\{1,\dots, n_2-1\}\nonumber.
\end{align}

Moreover, denoting the $i$-th principal curvature function of
$f_s^1$ (resp. $f_{\phi(s)}^2$) by $\lambda_i^s$ (resp. $\mu_i^s$),
the principal curvatures of $\Sigma$
are the functions
\begin{equation}\label{eq-principalcurvatures}
  k_{1}=\frac{\phi''}{\mathcal W^3}, \quad
  k_{1+i_1}=\frac{-\phi'}{\mathcal W}\lambda_{i_1}^s, \quad
  k_{n_1+i_2}=\frac{1}{\mathcal W}\mu_{i_2}^s,
\end{equation}
where $1\le i_1\le n_1-1$ and $1\le i_2\le n_2-1$. In particular,
the mean curvature of $\Sigma$  at $p=(p_1,p_2,s)$ is given by
\begin{equation} \label{eq-H}
H_\Sigma(p_1,p_2,s)=\frac{\phi''(s)}{\mathcal W^3(s)}-\frac{\phi'(s)}{\mathcal W(s)}H_{1}(p_1,s)+\frac{1}{\mathcal W(s)}H_2(p_2,\phi(s)),
\end{equation}
where $H_i(p_i,t)$ is the mean curvature of $f_t^i$ at $p_i$, $i=1,2$.

\begin{remark}  \label{rem:phi'constant}
It follows from~\eqref{eq-normal} that the angle function $\theta$ of
an $(f_s^1,f_{\phi(s)}^2)$-hyper\-surface of $\qq$ is constant if and only if
$\phi'$ is constant. If so, $\theta^2\ne1$. It should also be noticed that
the first principal curvature $k_{1}={\phi''}/{\mathcal W^3}$ of such a hypersurface depends only on $s$,
and coincides with the curvature function of the plane curve
$s\in I\mapsto (s,\phi(s))\in\R^2$.
\end{remark}

\begin{lemma} \label{lem-Hparallels}
Let $\Sigma\subset\qq$ be
an $(f_s^1,f_{\phi(s)}^2)$-hypersurface
of constant mean curvature. Then,
the families $f_s^1$ and $f_{\phi(s)}^2$  are both isoparametric.
\end{lemma}

\begin{proof}
Suppose that $H_\Sigma$ is constant, i.e., that it is independent of
$p_i$ and $s$. Fixing $s$ and $i\in\{1,2\}$, we conclude from~\eqref{eq-H}
that $H_j(p_j,s)$, $j\ne i$, is independent of $p_j$, which implies that
both families $f_s^1$ and $f_{\phi(s)}^2$ are isoparametric.
\end{proof}

It will be convenient to introduce the following concepts.

\begin{definition} \label{def-split}
Let $\Sigma\subset\qq$ be an oriented  hypersurface.
We say that the principal directions of $\Sigma$ \emph{split properly} at $p\in\Sigma$ if
$T(p)\ne0$, and there exists an orthonormal basis $$\{Z^1,\dots, Z^{n_1+n_2-1}\}\subset T_p\Sigma,$$
consisting of principal directions of $\Sigma$ at $p$, which can be
indexed so that:
\begin{itemize}[parsep=1ex]
\item[(i)] $Z^1=T/\|T\|$;
\item[(ii)] $Z^{i}\in\mathfrak{X}(\q_{\epsilon_1}^{n_1})\oplus\{0\} \,\,\, \forall i\in\{2,\dots, n_1\}$;
\item[(iii)] $Z^{i}\in\{0\}\oplus\mathfrak{X}(\q_{\epsilon_2}^{n_2}) \,\,\, \forall i\in\{n_1+1,\dots, n_1+n_2-1\}$.
\end{itemize}
In particular,
$$PZ^i=\pm Z^i \qquad \forall i\in\{2,\dots, n_1+n_2-1\}.$$
\end{definition}

\begin{definition} \label{def:distinguished}
Let $\Sigma$ be a hypersurface of $\qq$ with angle function
$\theta$. We say that a point $p\in\Sigma$
is \emph{distinguished}
if it satisfies either of the following conditions:
\begin{itemize}[parsep=1ex]
\item $\theta^2(p)=1$;
\item $\theta^2(p)<1$ and the principal directions of $\Sigma$ split properly at $p$.
\end{itemize}
\end{definition}

\begin{example}
It follows from~\eqref{eq-principaldirections} that  the principal directions of any
$(f_s^1,f_{\phi(s)}^2)$-hyper\-surface $\Sigma$ of $\qq$ split properly everywhere. In particular,
\emph{every point of $\Sigma$ is distinguished.}
\end{example}

\begin{example} \label{examp-horospherical}
Let $\Sigma$ be an $(f_s^1,f_{\phi(s)}^2)$-hypersurface of $\qq$,
where $\epsilon_i\le~0$.
Suppose that $\phi(s)=as, \,a>0$, $s\in\R$,  and that
$f_s^i$ is either a family of parallel hyperplanes
of $\R^{n_i}$ or of parallel horospheres of the hyperbolic space
$\h_{\epsilon_i}^{n_i}$.
Since $s$ varies in $\R$ and $\phi$ is linear, we have that $\Sigma$ is complete
and has constant angle $\theta\ne\pm1$. Moreover,
$f_s^i$ and $f_{s_*}^i$ are congruent for all $s, s_*\in\R$,
and  $f_s^i$ is invariant by a subgroup $\Gamma_i\subset {\rm iso}(\q_{\epsilon_i}^{n})$ of translations
of $\q_{\epsilon_i}^{n}$ for all $s\in\R$. Therefore, $\Sigma$
is homogeneous. In particular, it is isoparametric and has constant principal curvatures
(this last property also comes from~\eqref{eq-principalcurvatures}). We shall call such
a $\Sigma$ \emph{flat-horospherical} if some $\epsilon_i$ vanishes. Otherwise, it will be called
\emph{bi-horospherical}.

Note that, in the bi-horospherical case, if
$\phi(s)=s$, we have from~\eqref{eq-normal}
that the angle function of
$\Sigma$ vanishes identically.
If, in addition,  $n_1=n_2$ and $\epsilon_1=\epsilon_2<0$, then~\eqref{eq-H} implies
that $\Sigma$ is a minimal hypersurface of $\h_{\epsilon_1}^{n_1}\times\h_{\epsilon_2}^{n_2}$.
\end{example}

Next,  we characterize the hypersurfaces of
\,$\qq$ with constant angle $\theta=~\pm 1$.

\begin{proposition} \label{prop-Ni=0}
Let $\Sigma\subset\qq$ be an oriented hypersurface
with constant angle function $\theta=\pm 1.$
Then, locally, $\Sigma$ is an open set of one of the following
types of hypersurfaces:
\begin{itemize}[parsep=1ex]
\item[\rm (i)] $\Sigma_1^{n_1-1}\times\q_{\epsilon_2}^{n_2}$;
\item[\rm (ii)] $\q_{\epsilon_1}^{n_1}\times\Sigma_2^{n_2-1}$;
\end{itemize}
where $\Sigma_i^{n_i-1}$
is a  hypersurface of \,$\q_{\epsilon_i}^{n_i}$.
Furthermore, $\Sigma$ is of constant mean curvature (resp.~isoparametric, homogeneous) if and only if
$\Sigma_i^{n_i-1}$ is of constant mean curvature (resp.~isoparametric, homogeneous).
\end{proposition}

\begin{proof}
We can assume $\theta=1$, since the complementary case $\theta=-1$ can be treated analogously.
Then the unit normal $N$ of $\Sigma$ satisfies
$$N=N_1+N_2=N_1.$$

Set $\{X^1,\dots, X^{n_1-1}\}$ for a local orthonormal
frame  in $\{N_1\}^\perp\subset\mathfrak X(\q_{\epsilon_1}^{n_1})\oplus\{0\}$, and
$\{Y^1,\dots, Y^{n_2}\}$ for a local orthonormal
frame  in $\{0\}\oplus\mathfrak X(\q_{\epsilon_2}^{n_2})$. Denote by
$\mathfrak D_1$ and $\mathfrak D_2$ the
distributions $\{X^i\}$ and $\{Y^i\}$, respectively, and note that
$$\mathfrak X(\Sigma)=\mathfrak D_1\oplus\mathfrak D_2.$$

For all $i,j\in\{1,\dots, n_1-1\}$ and $k\in\{1,\dots, n_2\}$, one has
\[
\langle\nabla_{X^i}X^j,Y^k\rangle=\langle\overbar\nabla_{X^i}X^j,Y^k\rangle=0,
\]
so that $\mathfrak D_1$ is integrable.

Regarding $\mathfrak{D}_2$,
it is clear that $P|_{\mathfrak D_2}=-{\rm Id}_{\mathfrak D_2}$.
Since $PN=N_1=N$ and $P$ is parallel, we also have  $PA=A$.
Hence $A|_{\mathfrak D_2}=0$, i.e.,
$\mathfrak D_2$ is totally geodesic. Besides, for
any $i,j\in\{1,\dots, n_2\}$ and $k\in\{1,\dots, n_1-1\}$,
\[
\langle\nabla_{Y^i}Y^j,X^k\rangle=\langle\overbar\nabla_{Y^i}Y^j,X^k\rangle=0,
\]
which implies that $\mathfrak D_2$ is integrable and determines a totally geodesic
foliation of $\Sigma$, being each leaf locally isometric to $\q_{\epsilon_2}^{n_2}$.
Therefore $\Sigma$ is an open set of a hypersurface of type (i), as we wished to prove.
\end{proof}

\begin{example} \label{exam-tube}
Given $\epsilon\in\{-1,1\}$, set $I_\epsilon=(-\infty,-1)$ for
$\epsilon=-1$ and $I_\epsilon=(-1,1)$ for
$\epsilon=1.$
Choose  $t\in I_\epsilon$ and let $\Sigma_t$ be the tube
of radius $\arccos_\epsilon(\epsilon t/\sqrt 2)$ over the diagonal
$\{(p,p)\in\q_\epsilon^n\times\q_\epsilon^n\}$ of
$\q_\epsilon^n\times\q_\epsilon^n$, where
$\cosh_\epsilon=\cosh$ for $\epsilon=-1$ and $\cos_\epsilon=\cos$
for $\epsilon=1$. It can be proved that (see~\cite{gao,urbano})
\[
\Sigma_t=\{(p,q)\in\q_\epsilon^n\times\q_\epsilon^n\,;\, \la p,q\ra=t\},
\]
and that $\Sigma_t$ is homogeneous with  constant angle function $\theta=0$. Moreover,
at any point, the  principal directions of $\Sigma_t$ fail to split properly. In particular,
\emph{$\Sigma_t$ has no distinguished points.}
\end{example}

\begin{example} \label{exam-cui}
In~\cite{cui}, Cui obtained a one-parameter family of isoparametric hypersurfaces
in $\s_{\epsilon_1}^n \times \s_{\epsilon_2}^n$, $\epsilon_1,\epsilon_2>0$,
which generalizes Urbano's family $\Sigma_t$ from the above example.
As in Urbano's case, the hypersurfaces in Cui's family have constant angle equal to zero
and possess no distinguished points.
\end{example}


Let $\Sigma\subset\qq$ be a hypersurface with angle function $\theta\ne\pm1$.
Then
$$T=T_1+T_2=(1-\theta)N_1-(1+\theta)N_2\ne 0.$$

Since $T_i$ and $N_i$ are linearly dependent,
we can choose a local  orthonormal frame
$\mathscr F=\{U^1,\dots, U^{n-1}\}\subset\mathfrak X(\Sigma)$, $n=n_1+n_2$, such that
\[
U^1=\frac{T}{\|T\|}, \quad U^2, \dots ,U^{n_1}\in\mathfrak X(\q_{\epsilon_1}^{n_1})\oplus\{0\},
\quad U^{n_1+1}, \dots ,U^{n-1}\in\mathfrak \{0\}\oplus\mathfrak X(\q_{\epsilon_2}^{n_2}).
\]

By abuse of notation, we shall denote by $A=(a_{ij})$ the matrix of the shape operator
$A$ of $\Sigma$ with respect to $\{U^1,\dots, U^{n-1}\}$, that is,
\begin{equation} \label{eq-AUj}
AU^j=\sum_{i=1}^{n-1}a_{ij}U^i.
\end{equation}

If we assume further that $\theta$ is constant on $\Sigma$, we have from Lemma~\ref{lem-T&theta} that
\[
\overbar\nabla_{U^1}N=-AU^1=\frac{1}{2\|T\|}\nabla\theta=0,
\]
and also that
\[
\overbar\nabla_{U^1}U^1=\nabla_{U^1}U^1+\la AU^1,U^1\ra N=0.
\]

Hence we can assume
$\mathscr F$  parallel along the trajectories of $U^1$, i.e.,
\begin{equation} \label{eq-nablaU1=0}
\nabla_{U^1} U^i=\overbar\nabla_{U^1} U^i=0 \,\,\, \forall i\in\{1,\dots, n-1\}.
\end{equation}

For the next proposition, we decompose the matrix $A$ as follows.
Since $AU^1=0$, one has $a_{1k}=a_{k1}=0$ for all $k\in\{1,\dots, n-1\}$.
Then,  we can exclude the first row and the first
column of $A$. The resulting submatrix decomposes in blocks as
\begin{equation*}
\left[\begin{array}{c|c}
\mathcal A_1 & \mathcal A_2\\[1ex]
\hline
\\[-1.5ex]
\mathcal A_2^\top & \mathcal A_3
\end{array}\right],
\end{equation*}
where $\mathcal A_1$, $\mathcal A_2$  and $\mathcal A_3$ are the following
submatrices of $A$:
\begin{itemize}[parsep=1ex]
\item $\mathcal A_1=(a_{ij}), \,\, 2\le i,j\le n_1$;
\item $\mathcal A_2=(a_{ij}), \,\, 2\le i\le n_1, \, n_1+1\le j\le n-1$;
\item $\mathcal A_3=(a_{ij}), \,\, n_1+1\le i,j\le n-1$.
\end{itemize}

We also set
\begin{equation} \label{eq:Ttilde}
 c=\frac{1+\Theta}{1-\Theta}>0,
 \qquad
 \widetilde T=\frac{1}{1-\Theta}T,
 \qquad\text{and}\qquad
 \tau_\ell=-\epsilon_\ell\|N_\ell\|^2, \,\,\, \ell\in\{1,2\}.
\end{equation}

\begin{proposition} \label{prop-split}
Let $\Sigma\subset\qq$ be a hypersurface with constant angle
function $\theta\ne\pm1$. With the above notation, the following hold:
\begin{equation}\label{eq:system}
\left\{
\begin{array}{lllllll}
\mathcal A_1^2 &-& c\mathcal A_2\mathcal A_2^\top&-&\tau_1 I_{n_1-1}&=&\widetilde T(\mathcal A_1),\\[1ex]
\mathcal A_2^\top\mathcal A_2&-& c \mathcal A_3^2&+&c\tau_2I_{n_2-1}&=&\widetilde T(\mathcal A_3),\\[1ex]
\mathcal A_1\mathcal A_2 &-&c\mathcal A_2\mathcal A_3&&&=&\widetilde T(\mathcal A_2),
\end{array}
\right.
\end{equation}
where $ I_{n_i-1}$ is the identity matrix of order $n_i-1$, and
$\widetilde T(\mathcal A_i)$ is the matrix whose entries are
the derivatives $\widetilde T(a_{ij})$ of the entries $a_{ij}$ of $\mathcal A_i$.
\end{proposition}

\begin{proof}
We have that
\[
U^1=\frac1{\|T\|}(PN-\theta N)
\]
and
$$\nabla_{U^i}U^1=\overbar\nabla_{U^i}U^1-\la AU^1,U^i\ra=\overbar\nabla_{U^i}U^1 \,\,\, \forall i\ne 1.$$
Thus
\begin{equation} \label{eq-nablaUiU1}
\nabla_{U^i}U^1=\frac{1}{\|T\|}(-PAU^i+\theta AU^i) \,\,\, \forall i\in\{1,\dots, n-1\}.
\end{equation}

It follows from~\eqref{eq-nablaU1=0} and the Codazzi equation~\eqref{eq-codazzi}
for $X=U^1$ and $Y=U^i$, $i\in\{2,\dots, n-1\}$, that
\begin{equation} \label{eq-codazzi02}
\nabla_{U^1}AU^i+A\nabla_{U^i}U^1=(\overbar R(U^1,U^i)N)^\top=
(\overbar R_{\epsilon_1}(U^1,U^i)N+\overbar R_{\epsilon_2}(U^1,U^i)N)^\top.
\end{equation}

Choosing $i\in\{2,\dots, n_1\}$ and noticing that
\[
\overbar R(U^1,U^i)N=\overbar R_{\epsilon_1}(U^1,U^i)N=\epsilon_1\langle U_1^1,N_1\rangle U^i=
\frac{\epsilon_1}{\|T\|}\la T_1,N_1\rangle U^i=\frac{\epsilon_1\|T\|}{2}U^i,
\]
we have from~\eqref{eq-nablaU1=0}-\eqref{eq-codazzi02}  that
\begin{eqnarray}
U^1(a_{ii}) &=&\left\la\nabla_{U^1} AU^i,U^i\right\ra=-\la A\nabla_{U^i}U^1,U^i\ra+\frac{\epsilon_1\|T\|}{2} \nonumber\\
            &=&\frac1{\|T\|}\langle PAU^i-\theta AU^i,AU^i\rangle +\frac{\epsilon_1\|T\|}{2}\,\cdot \label{eq-a22}
\end{eqnarray}

Considering~\eqref{eq-AUj}, we have
\[
PAU^i=\sum_{k=2}^{n_1}a_{ki}U^k-\sum_{k=n_1+1}^{n-1}a_{ki}U^k,
\]
which yields
\[
\la PAU^i,AU^i\ra=\sum_{k=2}^{n_1}a_{ki}^2-\sum_{k=n_1+1}^{n-1}a_{ki}^2.
\]

Since $\|T\|^2=1-\theta^2$,  we conclude from~\eqref{eq-a22} that
\begin{equation} \label{eq-Taii}
(1-\theta)\sum_{k=2}^{n_1}a_{ki}^2-(1+\theta)\sum_{k=n_1+1}^{n-1}a_{ki}^2+\frac{\epsilon_1}2(1-\theta^2)= T(a_{ii})\,\,\, \forall i\in\{2,\dots, n_1\}.
\end{equation}

Analogously, for all $j\in\{n_1+1,\dots, n-1\}$, one has
\[
\overbar R(U^1,U^j)N=\overbar R_{\epsilon_2}(U_2^1,U^j)N_2=\epsilon_2\langle U_2^1,N_2\rangle U^j=
\frac{\epsilon_2}{\|T\|}\la T_2,N_2\rangle U^j=-\frac{\epsilon_2\|T\|}{2}U^j,
\]
and then, a calculation as the one above gives
\begin{equation} \label{eq-Tajj}
(1-\theta)\sum_{k=2}^{n_1}a_{kj}^2-(1+\theta)\sum_{k=n_1+1}^{n-1}a_{kj}^2-\frac{\epsilon_2}2(1-\theta^2)=T(a_{jj}) \,\,\, \forall j\in\{n_1+1,\dots, n-1\}.
\end{equation}

Finally,  taking the derivative $T(a_{ij})$ for $i\ne j$ and proceeding as above yields
\begin{equation} \label{eq-Taij}
(1-\theta)\sum_{k=2}^{n_1}a_{ik}a_{kj}-(1+\theta)\sum_{k=n_1+1}^{n-1}a_{ik}a_{kj}=T(a_{ij}) \,\,\, \forall i\ne j\in\{2,\dots, n-1\}.
\end{equation}

In matrix form, identities~\eqref{eq-Taii}--\eqref{eq-Taij} are the equations of the system
\[
\left\{
\begin{array}{lllllll}
(1-\theta)\mathcal A_1^2 &-&(1+\theta)\mathcal A_2\mathcal A_2^\top&+&\frac{\epsilon_1}{2}(1-\theta^2) I_{n_1-1}&=&T(\mathcal A_1),\\[1ex]
(1-\theta)\mathcal A_2^\top\mathcal A_2&-&(1+\theta)\mathcal A_3^2&-&\frac{\epsilon_2}{2}(1-\theta^2) I_{n_2-1}&=&T(\mathcal A_3),\\[1ex]
(1-\theta)\mathcal A_1\mathcal A_2 &-&(1+\theta)\mathcal A_2\mathcal A_3&&&=&T(\mathcal A_2).
\end{array}
\right.
\]
Dividing by $1-\theta$ and using~\eqref{eq-N1N2theta} and~\eqref{eq:Ttilde}
gives the system~\eqref{eq:system}.
\end{proof}


\section{Proof of Theorem~\ref{thm-isoparametric-->thetaconstant}} \label{sec-proofTh1}

In this section, we prove Theorem~\ref{thm-isoparametric-->thetaconstant}. First,
by means of Jacobi field theory, we determine, in terms of the geometry of the initial condition,
the shape operator
of the parallel hypersurfaces $\Sigma_r$ of a given isoparametric hypersurface $$\Sigma=\Sigma_0\subset\qq.$$

This will lead us to a real analytic function  $D=D(r)$ whose properties
will force the angle function $\theta$ of $\Sigma$ to be constant.
Since we are assuming $\Sigma$ connected, it will suffice to prove that $\theta$ is locally constant.
Therefore, after passing to an open set, we can and will assume that
\[
\theta^2<1\qquad \text{on}\,\,\, \Sigma.
\]

\subsection{The second fundamental form of the parallels}\label{sec41}
Define the map
\[
\begin{array}{cccc}
\Phi^r\colon & \Sigma &\to & \qq \\ [1ex]
             &     p  &\mapsto & \exp_p(rN_p),
\end{array}
\]
where $N_p$ is the unit normal field of $\Sigma$
at $p\in\Sigma$ and $\exp$ stands for
the exponential map of $\qq$. Again by the local nature of the result,
we can assume that, for a
small $\delta>0,$ and for all $r\in(-\delta,\delta),$
the map $\Phi^r$ is well defined and $\Sigma_r:=\Phi^r(\Sigma)$ is an embedded hypersurface
of $\qq$ which lies at distance $|r|$ from $\Sigma$.

Given $p\in\Sigma$,
let $\gamma_p(r)$ be the geodesic of $\qq$ such that
$\gamma_p(0)=p$ and $\gamma_p'(0)=N_p$, that is,
$\gamma_p(r)=\Phi^r(p)$. It is easily seen that
the unit normal to $\Sigma^r$ at $\gamma_p(r)$ is
$N(r):=\gamma_p'(r)$. In particular,
$N(r)$ is parallel along $\gamma_p$, and so
the angle function $\theta$ is independent of $r$. Indeed,
\[
N(\theta)=N\langle PN,N\rangle=\langle\overbar\nabla_NPN,N\rangle+\langle N,\overbar\nabla_NN\rangle=
\langle P\overbar\nabla_NN,N\rangle=0.
\]

Since $\theta^2<1$, we have that  $T=(1-\theta)N_1-(1+\theta)N_2$ never vanishes.
Moreover, from the considerations above,
$T$ is parallel along $\gamma_p$.
This, together with the fact that $T_i$ and $N_i$ are linearly dependent, allows
us to consider a parallel orthonormal frame $\{N,U^1,\dots, U^{n-1}\}$ along $\gamma$ such that
\[
U^1=\frac{T}{\|T\|}, \quad U^2, \dots ,U^{n_1}\in\mathfrak X(\q_{\epsilon_1}^{n_1})\oplus\{0\},
\quad U^{n_1+1}, \dots ,U^{n-1}\in\mathfrak \{0\}\oplus\mathfrak X(\q_{\epsilon_2}^{n_2}).
\]

For any $j\in\{1,\dots,n-1\}$, let $\zeta_j=\zeta_j(r)$ be the Jacobi field along $\gamma_p$ such
that
\begin{equation} \label{eq-initialconditions}
\zeta_j(0)=U^j(0) \quad\text{and}\quad \zeta'_j(0)=-AU^j(0),
\end{equation}
where $A$ is the shape operator of $\Sigma$.
Then, for any such $j$, $\zeta_j$ satisfies
\begin{equation} \label{eq-jacobi}
\zeta_j''+\overbar R(\gamma_p',\zeta_j)\gamma_p'=0.
\end{equation}
Besides, for all  $r\in(-\delta,\delta),$
one has
$\la\zeta_j(r),N(r)\ra=0$. Thus,
there exist smooth functions $b_{ij}=b_{ij}(r)$ such that
\begin{equation} \label{eq-zetacoordinates}
\zeta_j=\sum_{i=1}^{n-1}b_{ij}U^i, \,\,\, j\in\{1,\dots, n-1\}.
\end{equation}
Furthermore, since all $U^i$ are parallel along
$\gamma_p$, we have
\begin{equation}\label{eq-zeta''}
\zeta''_j=\sum_{i=1}^{n-1}b_{ij}''U^i, \,\,\, j\in\{1,\dots, n-1\}.
\end{equation}

Now, aiming the Jacobi equation~\eqref{eq-jacobi}, we compute
\begin{eqnarray}  \label{eq-jacobi02}
\overbar R(\gamma_p', \zeta_j)\gamma_p' &=&R_{\epsilon_1}\left(N_1,b_{1j}\frac{T_1}{\|T\|}+\sum_{i=2}^{n_1}b_{ij}U^i\right)N_1\\
                                        && \phantom{111111111111} +
                                        R_{\epsilon_2}\left(N_2,b_{1j}\frac{T_2}{\|T\|}+\sum_{i=n_1+1}^{n-1}b_{ij}U^i\right)N_2.\nonumber
\end{eqnarray}

Considering that
\[
\left\langle\frac{T_1}{\|T\|},N_1\right\rangle N_1=\frac{1-\theta}{\|T\|}\|N_1\|^2N_1=\|N_1\|^2\frac{T_1}{\|T\|},
\]
we have for the first summand in~\eqref{eq-jacobi02} that
\begin{eqnarray*}
R_{\epsilon_1}\left(N_1,b_{1j}\frac{T_1}{\|T\|}+\sum_{i=2}^{n_1}b_{ij}U^i\right)N_1 &=&
  \epsilon_1\|N_1\|^2\left(b_{1j}\frac{T_1}{\|T\|} +\sum_{i=2}^{n_1}b_{ij}U^i\right)\\
 && - \epsilon_1\left\langle b_{1j}\frac{T_1}{\|T\|},N_1\right\rangle N_1\\
 &=& \epsilon_1\|N_1\|^2\sum_{i=2}^{n_1}b_{ij}U^i.
\end{eqnarray*}

Analogously, for the second summand, we have
\[
R_{\epsilon_2}\left(N_2,b_{1j}\frac{T_2}{\|T\|}+\sum_{i=n_1+1}^{n-1}b_{ij}U^i\right)N_2=\epsilon_2\|N_2\|^2\sum_{i=n_1+1}^{n-1}b_{ij}U^i.
\]

Combining these two last equalities with~\eqref{eq-jacobi02} yields
\begin{equation}  \label{eq-jacobi03}
\overbar R(\gamma_p', \zeta_j)\gamma_p'=\epsilon_1\|N_1\|^2\sum_{i=2}^{n_1}b_{ij}U^i+\epsilon_2\|N_2\|^2\sum_{i=n_1+1}^{n-1}b_{ij}U^i.
\end{equation}
Recalling  the functions
$$\tau_\ell:=-\epsilon_\ell\|N_\ell\|^2, \qquad \ell\in\{1,2\},$$
introduced in~\eqref{eq:Ttilde},
equalities~\eqref{eq-jacobi} and~\eqref{eq-jacobi03}  give
\begin{equation} \label{eq-zeta''02}
\zeta_j''=\tau_1\sum_{i=2}^{n_1}b_{ij}U^i+\tau_2\sum_{i=n_1+1}^{n-1}b_{ij}U^i,  \,\,\, j\in\{1,\dots, n-1\}.
\end{equation}

Now, set $(a_{ij})$ for the (symmetric) matrix  $A$ of the shape operator of $\Sigma$ with respect to
the orthonormal basis $\{U^1(0),\dots, U^{n-1}(0)\}$, that is,
$$AU^j(0)=\sum_{i=1}^{n-1}a_{ij}U^i(0), \,\,\, j\in\{1,\dots, n-1\}.$$

Considering this last equality and comparing~\eqref{eq-zeta''} with~\eqref{eq-zeta''02}, we conclude
that $\zeta_j$ is a Jacobi field satisfying the initial conditions~\eqref{eq-initialconditions}
if and only if the coefficients
$b_{ij}$ are solutions of an initial value problem. Namely,
\begin{equation}\label{eq-IVP}
\left\{\begin{array}{rcll}
b_{1j}''(r)&=&0\,
;\\[1ex]
b_{ij}''(r)&=&\tau_1b_{ij}(r),&\text{if}\  \ 2\le i\le n_1\,;
\\[1ex]
b_{ij}''(r)&=&\tau_2b_{ij}(r),&\text{if}\ \ n_1+1\le i\le n-1\,;
\\[1ex]
b_{ij}(0)&=&\delta_{ij}
;\\[1ex]
b_{ij}'(0)&=&-a_{ij}.
\end{array}\right.
\end{equation}

It is easily checked that, if $\epsilon_1\epsilon_2\ne 0$,  the solutions of~\eqref{eq-IVP} are
\begin{equation}\label{eq-solutions}
\left\{\begin{array}{rcll}
b_{1j}(r)&=&\delta_{1j}-a_{1j}r;\\[1ex]
b_{ij}(r)&=&\delta_{ij}c_{\tau_1}(r)-a_{ij}s_{\tau_1}(r)&\text{if}\ \ \ 2\le i\le n_1\,;\\[1ex]
b_{ij}(r)&=&\delta_{ij}c_{\tau_2}(r)-a_{ij}s_{\tau_2}(r)&\text{if}\ \ n_1+1\le i\le n-1\,;
\end{array}\right.
\end{equation}
where, for $\ell\in\{1,2\}$, $s_{\tau_\ell}$ and $c_{\tau_\ell}$ are the functions:
\begin{equation}\label{def-sc}
s_{\tau_\ell}(r):=\left\{
\begin{array}{lcc}
\frac{1}{\sqrt{\tau_\ell}}\sinh(\sqrt\tau_\ell\, r) & \text{if} & \tau_\ell>0,\\[1ex]
\frac{1}{\sqrt{-\tau_\ell}}\sin(\sqrt{-\tau_\ell}\, r)  & \text{if} & \tau_\ell<0,
\end{array}
\right.
\quad
c_{\tau_\ell}:=\left\{
\begin{array}{lcc}
\cosh(\sqrt\tau_\ell\, r) & \text{if} & \tau_\ell>0,\\[1ex]
\cos(\sqrt{-\tau_\ell}\, r)  & \text{if} & \tau_\ell<0.
\end{array}
\right.
\end{equation}

Notice that the derivatives of  $s_{\tau_\ell}$ and $c_{\tau_\ell}$
satisfy:
\begin{equation} \label{eq-s'&c'}
s_{\tau_\ell}'(r)=c_{\tau_\ell}(r) \quad\text{and}\quad c_{\tau_\ell}'(r)=\tau_\ell s_{\tau_\ell}(r) \,\,\, \forall r\in\R.
\end{equation}

If $\epsilon_1=0\ne\epsilon_2$, the solutions of~\eqref{eq-IVP} are
\begin{equation}\label{eq-solutions02}
\left\{\begin{array}{rcll}
b_{ij}(r)&=&\delta_{ij}-a_{ij}r &\text{if}\ \ 2\le i\le n_1\,;\\[1ex]
b_{ij}(r)&=&\delta_{ij}c_{\tau_2}(r)-a_{ij}s_{\tau_2}(r)&\text{if}\  n_1+1\le i\le n-1.
\end{array}\right.
\end{equation}

Given $r\in(-\delta,\delta),$ let
$B(p,r)$  be the linear operator of $T_{\gamma_p(r)}\Sigma_r$ which takes
the basis $\mathcal B:=\{U^1(r), \dots, U^{n-1}(r)\}$ to
$\{\zeta_1(r), \dots, \zeta_{n-1}(r)\}$.
Considering~\eqref{eq-zetacoordinates} and the fact that each
$U^i$ is parallel along $\gamma_p$, we deduce that the matrix
of $B(p,r)$ (which we also denote by $B(p,r)$)
with respect to $\mathcal B$ is
$$B(p,r)=(b_{ij}(r)) \,\,\,  i,j\in\{1,\dots, n-1\},$$
where $b_{ij}$ are the functions defined in~\eqref{eq-solutions}.

In the above setting, Jacobi field theory applies, giving that
$B(p,r)$ is nonsingular for each $r\in(-\delta,\delta)$, and that the
shape operator  of $\Sigma^r$ at $p$ is  (see~\cite[Thm. 10.2.1]{olmos})
\begin{equation} \label{eq:A(p,r)}
A(p,r)=-B'(p,r)B(p,r)^{-1}.
\end{equation}

Since $\Sigma$ is isoparametric, the mean curvature $H(r)$ of $\Sigma^r$  is constant, i.e.,
it is independent of $p$.  Moreover, setting $$D(p,r)=\det B(p,r),$$
we have from~\eqref{eq:A(p,r)} that
\begin{equation}\label{eqH(r)}
H(r)={\rm tr}\,A(p,r)=-{\rm tr}\,(B'(p,r)B(p,r)^{-1})=-\frac{D'(p,r)}{D(p,r)}\cdot
\end{equation}
Since $D(0)=1$, it follows that
\[
D(p,r)=e^{-\int_{0}^{r}H(s)ds}.
\]
In particular, $D(p,r)$ is independent of $p$. So we can write
\[
D(r)=D(p,r)=\det B(p,r).
\]

Setting, for $\tau_1=0$, $$c_{\tau_1}(r)=1, \qquad s_{\tau_1}(r)=r,$$ and
considering the equalities~\eqref{eq-solutions}
and~\eqref{eq-solutions02} that define the functions $b_{ij}$,
we can write the matrix $B(p,r)=(b_{ij}(r))$  as

\begin{equation}\label{eq-B(r)}
B(r)=\small{
\left[
\begin{array}{c|c}
1-a_{11}r& \begin{array}{ccc}
-a_{12}r&\dots&-a_{1(n-1)}r
\end{array}\\
\hline
\begin{array}{c}
-a_{12}s_{\tau_1}(r)^{\phantom {2^1}}\\
\vdots\\
-a_{1n_1}s_{\tau_1}(r)_{\phantom{x_1}}
\end{array}& \delta_{ij}c_{\tau_1}(r)-a_{ij}s_{\tau_1}(r)\\
\hline
\begin{array}{c}
-a_{1(n_1+1)}s_{\tau_2}(r)^{\phantom {2^1}}\\
\vdots\\
-a_{1(n-1)}s_{\tau_2}(r)_{\phantom{x_1}}
\end{array}& \delta_{ij}c_{\tau_2}(r)-a_{ij}s_{\tau_2}(r)
\end{array}
\right]}.
\end{equation}
\smallskip

A direct computation shows that, given
$$u\in\{0,\dots, n_1\} \qquad\text{and}\qquad v\in\{0,\dots, n_2-1\},$$
there exist functions
$$\alpha_{u,v,p}=\alpha_{u,v,p}(a_{ij}) \quad\text{and}\quad \beta_{u,v,p}=\beta_{u,v,p}(a_{ij}), \quad  p\in\Sigma,  \quad A=(a_{ij}),$$
such that
\begin{equation} \label{eq-D&cs}
D(r)=\left\{
\begin{array}{ll}
\displaystyle\sum_{u=0}^{n_1-1}\sum_{v=0}^{n_2-1}\left(\alpha_{u,v}+
r\beta_{u,v}\right)\cone^{n_1-1-u}(r)\sone^u(r)\ctwo^{n_2-1-v}(r)\stwo^v(r), & \tau_1\tau_2\ne 0,\\[3.5ex]
\displaystyle\sum_{u=0}^{n_1}\sum_{v=0}^{n_2-1}\alpha_{u,v}r^u c_{\tau_2}^{n_2-1-v}(r)s_{\tau_2}^v(r), & \tau_1=0\ne\tau_2.
\end{array}
\right.
\end{equation}

For $\tau\ne 0$, define
\[
 \omega(\tau)=
 \begin{cases}
  \sqrt{\tau},&\tau>0,\\[1ex]
  \mathrm{i}\sqrt{-\tau},&\tau<0,
 \end{cases}
\]
and note that
\begin{equation} \label{eq:cs&exp}
 c_\tau(r)=\frac{e^{\omega(\tau)r}+e^{-\omega(\tau)r}}2,
 \qquad
 s_\tau(r)=\frac{e^{\omega(\tau)r}-e^{-\omega(\tau)r}}{2\omega(\tau)}\cdot
\end{equation}

Hence, if we define the finite sets
\begin{equation}\label{def-N_i}
\left\{
\begin{array}{ll}
\mathcal N_\ell:=\{n_\ell-1-2u \,;\, u=0,\dots, n_\ell-1\}, & \ell\in\{1,2\},\\[1.5ex]
E(p):=
 \left\{a_1\sqrt{\tau_1(p)}
       +a_2\sqrt{\tau_2(p)} \,;\,
       a_\ell\in\mathcal N_\ell\right\}, &  p\in\Sigma,
 \end{array}
\right.
\end{equation}
and substitute~\eqref{eq:cs&exp} into~\eqref{eq-D&cs}, we get
\begin{equation}\label{eq:ca-expansion}
 D(r)=\sum_{\lambda\in E(p)}\varphi_{\lambda,p}(r)e^{\lambda r},
 \qquad \varphi_{\lambda,p}\in\mathbb C[r],
\end{equation}
where the polynomial coefficients $\varphi_{\lambda,p}$ are affine when $\tau_1\tau_2\ne0$,
and have degree at most $n_1$ when $\tau_1=0$.

We shall use the standard fact that, if $E\subset\mathbb C$ is finite and its
elements are pairwise distinct, then, for $p_\lambda\in\mathbb C[r]$,
\begin{equation} \label{eq:elementary}
\sum_{\lambda\in E}p_\lambda(r)e^{\lambda r}=0
\qquad\Longrightarrow\qquad p_\lambda=0 \quad \forall \lambda\in E.
\end{equation}
This follows by applying
successively suitable powers of the differential operators
$\frac{d}{dr}-\lambda$ to this identity.

Because $D$ is independent of the base point, taking $p\neq q\in\Sigma$, and
setting $\varphi_{\lambda,p}=0$ when $\lambda\notin E(p)$, and similarly for $q$,
we have
\[
0=D(r)-D(r)
=\sum_{\lambda\in E(p)\cup E(q)}
\bigl(\varphi_{\lambda,p}(r)-\varphi_{\lambda,q}(r)\bigr)e^{\lambda r}.
\]

 Since the indices $\lambda$ belong to a finite set and are pairwise distinct,~\eqref{eq:elementary}
 yields
\[
\varphi_{\lambda,p}=\varphi_{\lambda,q}
\qquad\forall\lambda\in E(p)\cup E(q).
\]
Therefore the reduced support
\begin{equation}\label{eq:sec4-38}
\Lambda:=\{\lambda\in E(p):\varphi_{\lambda,p}\not\equiv0\}
\end{equation}
is independent of $p$. We may consequently write
\begin{equation}\label{exp-D}
D(r)=\sum_{\lambda\in\Lambda}\varphi_\lambda(r)e^{\lambda r},
\end{equation}
where $\varphi_\lambda\in\mathbb C[r]$ is independent of $p$. Note that, since
$D(0)=1$, one has $\Lambda\neq\varnothing$.


\subsection{$\boldsymbol {D(r)}$ as a polynomial}
As we observed earlier, the function $D(r)$ will play a crucial role in the proof of
Theorem~\ref{thm-isoparametric-->thetaconstant}. Curiously, the central question about this
function is whether it can be a polynomial in $r$.
On this matter, we  note that, if
$$D(r)=\psi(r)\in\mathbb C[r],$$
then~\eqref{exp-D} gives
\begin{equation} \label{eq:suum}
\sum_{\lambda\in\Lambda}\varphi_\lambda(r)e^{\lambda r}-\psi(r)=0.
\end{equation}

Suppose that $0\not\in\Lambda$. In this case,
the exponents in $\lambda\in\Lambda\cup\{0\}$
are pairwise distinct.  Hence, applying~\eqref{eq:elementary} to~\eqref{eq:suum} gives
\[
\varphi_\lambda\equiv \psi\equiv 0 \qquad \forall \lambda\in\Lambda,
\]
contradicting $D\not\equiv 0$. Thus $0\in\Lambda$ and~\eqref{eq:suum}
becomes
\[
\sum_{\lambda\in\Lambda\setminus\{0\}}\varphi_\lambda(r)e^{\lambda r}+(\varphi_0-\psi)(r)=0.
\]
Again by~\eqref{eq:elementary}, we have
\[
\psi\equiv\varphi_0 \qquad\text{and}\qquad \varphi_\lambda\equiv 0 \quad \forall\lambda\in\Lambda\setminus\{0\},
\]
and so $\Lambda=\{0\}.$

 Conversely, if $\Lambda=\{0\}$, it follows from~\eqref{exp-D} that
 \[
 D\equiv\varphi_0.
 \]

 Therefore, we have the following result.

 \begin{lemma} \label{lem:polynomial-equivalence}
 For $r$ near zero, the function $D=D(r)$ is a polynomial in $r$ if and only if $\Lambda=\{0\}$.
\end{lemma}

\begin{example}\label{minimal-bi-horo}
The minimal {bi-horospherical} hypersurfaces of Example~\ref{examp-horospherical}
realize $D(r)$ as a polynomial.
Take, for simplicity, $\epsilon_1=\epsilon_2=-1$, but $n_1$ and $n_2$ arbitrary.
Applying~\eqref{eq-normal} and~\eqref{eq-principalcurvatures} to such a bi-horospherical
hypersurface, we have
$$\tau_1=\frac{a^2}{1+a^2}, \qquad \tau_2=\frac{1}{1+a^2},$$
and its second fundamental form is
$$
A=\left[\begin{array}{ccc}
0&0&0\\
0&-\sqrt{\tau_1}I_{n_1-1}&0\\
0&0&\sqrt{\tau_2}I_{n_2-1}
\end{array}\right].
$$
Therefore
\begin{eqnarray*}
D(r)&=&\left[\cosh(r\sqrt{\tau_1})+\sinh(r\sqrt{\tau_1})\right]^{n_1-1}\left[\cosh(r\sqrt{\tau_2})-\sinh(r\sqrt{\tau_2})\right]^{n_2-1}\\
&=&e^{\left[(n_1-1)\sqrt{\tau_1}-(n_2-1)\sqrt{\tau_2}\right]r}.
\end{eqnarray*}

In particular, taking $a=(n_2-1)/(n_1-1)$ gives $D\equiv 1$.
Note that, in this case, the hypersurface is minimal, and the pair
$(\sqrt{\tau_1},\sqrt{\tau_2})$ is $(n_1,n_2)$-dependent, according to the next definition.
\end{example}

\begin{definition}\label{n1n2dep}
Given integers $n_1, n_2\ge 2$, we shall say that $(z_1,z_2)\in\mathbb C^2$ is
$(n_1,n_2)$-\emph{independent} if, for
$$a_\ell, b_\ell\in\mathcal N_\ell=\{n_\ell-1-2u \,;\, u=0,\dots, n_\ell-1\}, \quad \ell\in\{1,2\},$$
the following implication holds:
\begin{equation}\label{eq-implications}
(a_1-b_1)z_1+(a_2-b_2)z_2=0 \qquad \Longrightarrow \qquad a_\ell=b_\ell, \quad \ell\in\{1,2\}.
\end{equation}
Otherwise, the  pair $(z_1,z_2)$ will be said to be $(n_1,n_2)$-\emph{dependent}.
\end{definition}

\begin{lemma}\label{lemma-n1n2dep}
Let  $n_1, n_2\ge 2$ be integers, and assume $(\tau_1,\tau_2)\in\mathbb R^2$ with $\tau_1\tau_2\neq 0$.
Then the following are equivalent:
\begin{itemize}[parsep=1ex]
\item[\rm (i)] $\big(\sqrt{\tau_1},\sqrt{\tau_2}\big)$ is $(n_1,n_2)$-independent,
\item[\rm (ii)] the map $(a_1,a_2)\in\mathcal N_1\times\mathcal N_2\mapsto a_1\sqrt{\tau_1}+a_2\sqrt{\tau_2}\in\mathbb C$ is injective,
\item[\rm (iii)] for every $a_{\ell}\in\{1,\dots,n_{\ell}-1\}$, one has $a_1^2\tau_1\neq a_2^2\tau_2$,
\item[\rm (iv)] for integers $u, v$ such that $0\leq u\leq n_1-1,$ and $0\leq v\leq n_2-1,$ the functions
$$
f_{u,v}(r)=c_{\tau_1}^{\,n_1-1-u}(r)s_{\tau_1}^{\,u}(r)
c_{\tau_2}^{\,n_2-1-v}(r)s_{\tau_2}^{\,v}(r)
$$
are linearly independent.
\end{itemize}
In particular,  if $\tau_1\tau_2<0$, then $(\sqrt{\tau_1},\sqrt{\tau_2})$ is $(n_1,n_2)$-independent.
\end{lemma}

\begin{proof}
{\rm(i)} $\Leftrightarrow$ {\rm(ii)}. It follows directly from
the definition of \((n_1,n_2)\)-independence.

\medskip
\noindent
{\rm(ii)} $\Leftrightarrow$ {\rm(iii)}.
Since
$\mathcal N_\ell=\{n_\ell-1-2j:j=0,\dots,n_\ell-1\},$
the map in {\rm(ii)} fails to be injective if and only if there exist
integers \(a_1,a_2\) satisfying
$0<|a_1|\leq n_1-1$, $0<|a_2|\leq n_2-1$, and
such that
\[
a_1\sqrt{\tau_1}+a_2\sqrt{\tau_2}=0.
\]
Indeed, the differences between two elements of \(\mathcal N_\ell\)
are precisely the even integers \(2a_\ell\) satisfying
\(|a_\ell|\leq n_\ell-1\). The preceding equality implies
\[
a_1^2\tau_1=a_2^2\tau_2.
\]

Conversely, if this last equality holds, then
\[
a_1\sqrt{\tau_1}=\pm a_2\sqrt{\tau_2},
\]
and, after possibly changing the sign of $a_2$, the map in
{\rm(ii)} is not injective. This shows that  {\rm(ii)} and {\rm(iii)} are
equivalent.

\medskip
\noindent
(ii) $\Leftrightarrow$ (iv).
For
$0\leq u\leq n_1-1$ and $0\leq v\leq n_2-1,$
set
\[
g_{u,v}(r)=e^{(n_1-1-2u)\sqrt{\tau_1}r+(n_2-1-2v)\sqrt{\tau_2}r}.
\]
Using
\[
e^{\sqrt{\tau_\ell}r}=c_{\tau_\ell}(r)+
\sqrt{\tau_\ell}\,s_{\tau_\ell}(r),
\qquad
e^{-\sqrt{\tau_\ell}r}=
c_{\tau_\ell}(r)-\sqrt{\tau_\ell}\,s_{\tau_\ell}(r),
\]
we obtain
\[
\begin{aligned}
g_{u,v}(r)
={}&
\bigl(c_{\tau_1}(r)+\sqrt{\tau_1}s_{\tau_1}(r)\bigr)^{n_1-1-u}
\bigl(c_{\tau_1}(r)-\sqrt{\tau_1}s_{\tau_1}(r)\bigr)^u
\\
&\times
\bigl(c_{\tau_2}(r)+\sqrt{\tau_2}s_{\tau_2}(r)\bigr)^{n_2-1-v}
\bigl(c_{\tau_2}(r)-\sqrt{\tau_2}s_{\tau_2}(r)\bigr)^v.
\end{aligned}
\]

Expanding these expressions using the binomial formula
gives that every  $g_{u,v}$ belongs to
the vector space spanned by the functions $f_{u,v}$.
Moreover, it is easily checked that for $m\geq 1$ and $a\neq 0$, the homogeneous polynomials
\[
p(x,y)=(x+ay)^{m-j}(x-ay)^j,
\qquad j=0,\dots,m,
\]
are linearly independent. Applying this observation with
\[
(m,a)=(n_1-1,\sqrt{\tau_1})
\quad\text{and}\quad
(m,a)=(n_2-1,\sqrt{\tau_2})
\]
shows that the families \(\{f_{u,v}\}\) and \(\{g_{u,v}\}\) are related
by an invertible change of basis. Consequently,
\[
\{f_{u,v}\}\ \text{is linearly independent}
\quad\Longleftrightarrow\quad
\{g_{u,v}\}\ \text{is linearly independent}.
\]

Finally, a finite family of exponential functions
\(\{e^{\lambda r}\}\) is linearly independent if and only if its
exponents are pairwise distinct. Therefore, \(\{g_{u,v}\}\) is
linearly independent if and only if the numbers
\[
(n_1-1-2u)\sqrt{\tau_1}
+
(n_2-1-2v)\sqrt{\tau_2}
\]
are pairwise distinct. This is precisely the injectivity condition in
{\rm(ii)}.

\medskip
If \(\tau_1\tau_2<0\), the equality
\(a_1^2\tau_1=a_2^2\tau_2\) is impossible for nonzero
\(a_1,a_2\). Hence {\rm(iii)} holds, and the final assertion follows.
\end{proof}

In view of Lemma~\ref{lemma-n1n2dep}, we introduce the following concept.

\begin{definition}\label{redundant}
Given integers $n_1, n_2\ge 2$, and $\epsilon_1\epsilon_2\neq 0$,
we define the \emph{redundant set} $\mathcal R$ of $\Sigma$ as
$$
\mathcal R:=\{p\in\Sigma \,;\, \big(\sqrt{\tau_1(p)},\sqrt{\tau_2(p)}\big)\ \text{is}\ (n_1,n_2)\text{-dependent}\}.
$$
If $\epsilon_1\epsilon_2=0$, we set $\mathcal R=\emptyset$.
\end{definition}

\begin{remark} \label{rem:Rclosed}
Lemma~\ref{lemma-n1n2dep} shows that $\mathcal R=\emptyset$
if $\epsilon_1\epsilon_2<0$, for, in this case, $\tau_1\tau_2<0$.
Since, by definition, $\mathcal R=\emptyset$ if $\epsilon_1\epsilon_2=0$, we then have
\[
\mathcal R\ne\emptyset \qquad\Longrightarrow\qquad  \epsilon_1\epsilon_2> 0.
\]
Moreover, since $\tau_1$ and $\tau_2$ are continuous, the equivalence (i) $\Leftrightarrow$ (iii) in
Lemma~\ref{lemma-n1n2dep} implies that $(n_1,n_2)$-independence is an open condition.
In particular, $\mathcal R$ is closed.
\end{remark}


\begin{lemma}\label{lem:parity}
Suppose that $D(r)$ is a polynomial and $\Sigma\setminus\mathcal R\neq\emptyset$.
Then
$$\epsilon_{\ell}\neq 0 \qquad\Longrightarrow\qquad n_{\ell} \,\, \text{odd}, \quad\ell\in\{1,2\}. $$
\end{lemma}

\begin{proof}
Let us first note that
\begin{equation} \label{eq-implic}
0\in\mathcal N_\ell \qquad\Longleftrightarrow\qquad n_\ell \,\, \text{is odd}.
\end{equation}

Since  $D(r)$ is a polynomial, it follows from Lemma~\ref{lem:polynomial-equivalence}
that $\Lambda=\{0\}$, and so $0\in E(p)$ for every $p\in\Sigma$.
Hence, for each $\ell\in\{1,2\}$, there exists $a_\ell\in\mathcal N_\ell$
such that
\begin{equation} \label{eq-implic002}
a_1\sqrt{\tau_1}+a_2\sqrt{\tau_2}=0.
\end{equation}

Fix $p\in\Sigma\setminus\mathcal R$. If
$\epsilon_1\epsilon_2\ne 0$, then
$\big(\sqrt{\tau_1(p)},\sqrt{\tau_2(p)}\big)$ is $(n_1,n_2)$-independent,
and it follows from Lemma~\ref{lemma-n1n2dep}(ii) that $a_1=a_2=0$. This and~\eqref{eq-implic} show that
both $n_1$ and $n_2$ are odd.

Now suppose that $\epsilon_1=0\ne\epsilon_2$. Then $\tau_1=0\ne\tau_2$, and it follows
from~\eqref{eq-implic002} that $a_2=0$. Thus, again from~\eqref{eq-implic}, we conclude
that $n_2$ is odd.
\end{proof}

Next, we characterize the symmetric matrices $A$ whose corresponding
function $D(r)$ is a polynomial.
Then, we show that such an $A$ fails to be the shape
operator of any hypersurface of \,$\qq$.

With respect to our frame $\{U^1,\dots, U^{n_1+n_2-1}\},$
it will be convenient to decompose the matrix $A$ in blocks as
follows (note that this differs slightly from the notation in Section~\ref{sec41}):
If $\epsilon_1\epsilon_2\neq0$, we set
\begin{equation}\label{Aneq0}
A=\left[\begin{array}{c|c|c}
\mathcal A_{00}&\mathcal A_{01}&\mathcal A_{02}\\[1ex]
\hline\\[-2ex]
\mathcal A_{01}^{\top}&\mathcal A_{11}&\mathcal A_{12}\\
[1ex]\hline\\[-2ex]
\mathcal A_{02}^{\top}&\mathcal A_{12}^{\top}&\mathcal A_{22}\\[1ex]
\end{array}\right],
\end{equation}
where, denoting by $\R^{k,\ell}$ the space of real $k\times\ell$ matrices,
\[
\mathcal A_{00}=[a_{11}]\in\mathbb R^{1,1}, \quad \mathcal A_{0,\ell}\in\mathbb R^{1,n_{\ell}-1}, \quad
\mathcal A_{\ell,\ell}\in\mathbb R^{n_{\ell}-1,n_{\ell}-1}, \quad \mathcal A_{12}\in\mathbb R^{n_1-1,n_2-1}.
\]

When $\epsilon_1=0\neq\epsilon_2$, we merge some blocks and write
\begin{equation}\label{A=0}
A=\left[\begin{array}{c|c}
\widetilde A_{00}&{\mathcal A}_{02}\\[1ex]
\hline\\[-2ex]
 {\mathcal A}_{02}^{\top}&\mathcal A_{22}\\[1ex]
\end{array}\right],
\end{equation}
where
\[
\widetilde A_{00}\in\mathbb R^{n_{1},n_1}, \quad \mathcal A_{22}\in\mathbb R^{n_{2}-1,n_{2}-1}, \quad
{\mathcal A}_{02}\in\mathbb R^{n_1,n_2-1}.
\]

\begin{lemma}\label{D-pol}
Suppose that $\Sigma\setminus\mathcal R\neq \emptyset$. If $D(r)$ is a polynomial,
the following hold:
\begin{itemize}[parsep=1ex]
\item[\rm (i)]  $\epsilon_1\epsilon_2\neq 0$ \,\, $\Rightarrow$ \,\, $\epsilon_1,\epsilon_2<0$ \, and
\begin{equation*}
\mathcal A_{\ell\ell}^2=\tau_\ell I_{n_\ell-1},\quad{\rm tr}\mathcal A_{\ell\ell}=0,\quad \mathcal A_{01}=\mathcal A_{02}=\mathcal A_{12}=0.
\end{equation*}
\item[\rm (ii)] $\epsilon_1=0\neq \epsilon_2$ \,\, $\Rightarrow$ \,\, $\epsilon_2<0$ \, and
\begin{equation*}
\mathcal A_{22}^2=\tau_2 I_{n_2-1},\quad{\rm tr}\,\mathcal A_{22}=0,\quad\mathcal A_{02}=0.
\end{equation*}
\end{itemize}
\end{lemma}
\begin{proof}
(i) Fix  $p\in\Sigma\setminus\mathcal R$
and set
$$k_{\ell}:=\frac{n_{\ell}-1}{2}, \qquad \ell\in\{1,2\}.$$
By Lemma~\ref{lem:parity}, every  $k_{\ell}$ is an integer.

Now we introduce the following auxiliary variables:
\begin{equation*}
z_{\ell}=z_\ell(r)=
 \begin{cases}
  e^{2r\sqrt{\tau_{\ell}(p)}},&\text{if}\ \tau_{\ell}>0,\\[1ex]
  e^{2\mathrm{i}r\sqrt{-\tau_{\ell}(p)}},&\text{if}\ \tau_{\ell}<0.
 \end{cases}
\end{equation*}

In this way,  the functions $c_{\tau_{\ell}}$ and $s_{\tau_{\ell}}$
introduced in~\eqref{def-sc} can be rewritten as
\begin{equation}\label{csz}
c_{\tau_{\ell}}(r)=\frac{z_{\ell}+1}{2\sqrt{z_{\ell}}},\quad
s_{\tau_{\ell}}(r)=\frac{z_{\ell}-1}{2\sqrt{\tau_{\ell}z_{\ell}}}\cdot
\end{equation}

We modify the matrix $B(r)$ defined in~\eqref{eq-B(r)} as follows: first substitute \eqref{csz},
then keep the first row unchanged, multiply all rows
from the second to $n_1$-th by $2z_1^{\frac 12}$,
and all the others by $2z_2^{\frac 12}$. We then get the following matrix:
\begin{equation}\label{deftildeB}
\tilde B(r,z_1,z_2):=\left[\begin{array}{c|c|c}
1-r\mathcal A_{00}&-r\mathcal A_{01}&-r\mathcal A_{02}\\[1ex]
\hline\\[-2ex]
-\frac{z_1-1}{\sqrt{\tau_1}}\mathcal A_{01}^{\top}&(z_1+1)I_{n_1-1}-
\frac{z_1-1}{\sqrt{\tau_1}}\mathcal A_{11}&-\frac{z_1-1}{\sqrt{\tau_1}}\mathcal A_{12}\\
[1ex]\hline\\[-2ex]
-\frac{z_2-1}{\sqrt{\tau_2}}\mathcal A_{02}^{\top}&-\frac{z_2-1}{\sqrt{\tau_2}}
\mathcal A_{12}^{\top}&(z_2+1)I_{n_2-1}-\frac{z_2-1}{\sqrt{\tau_2}}\mathcal A_{22}\\[1ex]
\end{array}\right]
\end{equation}

Let $\tilde D(r,z_1,z_2):=\det \tilde B(r,z_1,z_2)$. By construction, we have
\begin{equation}\label{DvsD}
\tilde D(r,z_1,z_2)=2^{n_1+n_2-2}z_1^{k_1}z_2^{k_2} D(r).
\end{equation}
Then, by the  multilinearity of the determinant,
$$
\tilde D(r,z_1,z_2)=\sum_{a,b}M_{ab}(r)z_1^az_2^b, \qquad M_{ab}\in \mathbb C[r],
$$
where $0\leq a\leq n_1-1$ and $0\leq b\leq n_2-1.$

The fact that $p\notin\mathcal R$, together with
Lemma~\ref{lemma-n1n2dep}, implies that the corresponding exponential
factors $z_1^az_2^b$ are pairwise distinct.
The implication~\eqref{eq:elementary} then gives
$$
D(r)\ \text{polynomial}\quad\Longleftrightarrow\quad M_{ab}\equiv 0\ \text{for any}\ (a,b)\neq(k_1,k_2).
$$
Therefore
\begin{equation}\label{eqD1}
\tilde D(r,z_1,z_2)=M_{k_1,k_2}(r)z_1^{k_1}z_2^{k_2}.
\end{equation}
Considering $\tilde B(r,1,1)$, we have that
\begin{equation}\label{eqD2}
M_{k_1,k_2}(r)=\tilde D(r,1,1)=2^{n_1+n_2-2}(1-a_{11}r).
\end{equation}

Now let us fix $z_1=1$, so we can focus on the blocks $\mathcal A_{ 02}$ and $\mathcal A_{22}$:
\begin{eqnarray}
\nonumber\tilde D(r,1,z_2)&=&\det\left[\begin{array}{c|c|c}
1-r\mathcal A_{00}&-r\mathcal A_{01}&-r\mathcal A_{02}\\[1ex]
\hline\\[-2ex]
0&2I_{n_1-1}&0\\
[1ex]\hline\\[-2ex]
-\frac{z_2-1}{\sqrt{\tau_2}}\mathcal A_{02}^{\top}&-\frac{z_2-1}{\sqrt{\tau_2}}\mathcal A_{12}^{\top}&(z_2+1)I_{n_2-1}-\frac{z_2-1}{\sqrt{\tau_2}}\mathcal A_{22}\\[1ex]
\end{array}\right]\\[1.5ex]
\nonumber&=&2^{n_1-1}\det\left[\begin{array}{c|c}
1-r\mathcal A_{00}&-r\mathcal A_{02}\\[1ex]
\hline\\[-2ex]
-\frac{z_2-1}{\sqrt{\tau_2}}\mathcal A_{02}^{\top}&(z_2+1)I_{n_2-1}-\frac{z_2-1}{\sqrt{\tau_2}}\mathcal A_{22}\\[1ex]
\end{array}\right]\\[1ex]
\nonumber&=&2^{n_1-1}(1-a_{11}r)\det\left[(z_2+1)I_{n_2-1}-\frac{z_2-1}{\sqrt{\tau_2}}\tilde{\mathcal A}_{22}\right]\\
\label{eqD3}&=&2^{n_1-1}(1-a_{11}r)(z_2+1)^{n_2-1}\det\left[I_{n_2-1}-\frac{z_2-1}{\sqrt{\tau_2}(z_2+1)}\tilde{\mathcal A}_{22}\right],
\end{eqnarray}
where
$$
\tilde{\mathcal A}_{22}=\mathcal A_{22}+\frac{r}{1-a_{11}r}\mathcal A_{02}^{\top}\mathcal A_{02}.
$$

In the calculation above, we used that $1-ra_{11}\neq 0$ for $|r|$ small. This allowed
us to apply,  in the third equality, the Schur complement
for the determinant of a block matrix. Comparing~\eqref{eqD1},~\eqref{eqD2} and~\eqref{eqD3}, we have
$$
\det\left[I_{n_2-1}-\frac{z_2-1}{\sqrt{\tau_2}(z_2+1)}\tilde{\mathcal A}_{22}\right]=\left(\frac{4z_2}{(z_2+1)^2}\right)^{k_2}.
$$

We introduce two further auxiliary variables:
$$\zeta_{\ell}=\frac{z_\ell-1}{\sqrt{\tau_\ell}(z_\ell+1)}, \qquad \ell\in\{1,2\}.$$
Note that $1-\tau_\ell\zeta_\ell^2={4z_\ell}/{(z_\ell+1)^2}$, hence
$$
\det\left[I_{n_2-1}-\zeta_2\tilde{\mathcal A}_{22}\right]=\left(1-\tau_2\zeta_2^2\right)^{k_2},
$$
which implies that the eigenvalues of $\tilde{\mathcal A}_{22}$ are $\pm\sqrt{\tau_2}$,
each one with multiplicity $k_2$. In particular its trace vanishes.
Thus, by the definition of $\tilde{\mathcal A}_{22}$, for any $r$ we have
$$
0=\text{tr}\tilde{\mathcal A}_{22}=\text{tr}\mathcal A_{22}+\frac{r}{1-a_{11}r}|\mathcal A_{02}|^2.
$$
Therefore $\text{tr}\mathcal A_{22}=0$ and $\mathcal A_{02}=0$. So $\tilde{\mathcal A}_{22}=\mathcal A_{22}$,
which yields $\mathcal A_{22}^2=\tau_2I_{n_2-1}$ and $\text{tr}\mathcal A_{22}=0$.
Moreover, since $\mathcal A_{22}$ is real symmetric, its eigenvalues must be real.
In our setting this means $\tau_2>0$, equivalently $\epsilon_2<0$.
Repeating these arguments for the first factor we obtain
\[
\epsilon_1<0, \quad \mathcal A_{11}^2=\tau_1 I_{n_1-1}, \quad
{\rm tr}\mathcal A_{11}=0, \quad  \mathcal A_{01}=0.
\]

To  conclude the proof of (i), let us  show that $\mathcal A_{12}=0$.
To that end, define
$$
B_*(z_1,z_2):=\left[\begin{array}{c|c}
(z_1+1)I_{n_1-1}-\frac{z_1-1}{\sqrt{\tau_1}}\mathcal A_{11}&-\frac{z_1-1}{\sqrt{\tau_1}}\mathcal A_{12}\\
[1ex]\hline\\[-2ex]
-\frac{z_2-1}{\sqrt{\tau_2}}\mathcal A_{12}^{\top}&(z_2+1)I_{n_2-1}-\frac{z_2-1}{\sqrt{\tau_2}}\mathcal A_{22}\\[1ex]
\end{array}\right]
$$
and
$$
D_*(z_1,z_2):=\det B_*(z_1,z_2).
$$

From the above computations,  we have $\tilde D(r,z_1,z_2)=(1-ra_{11})D_*(z_1,z_2)$. Hence
$$
D_*(z_1,z_2)=2^{n_1+n_2-2}z_1^{k_1}z_2^{k_2}.
$$

We can operate on $B_*$ as before: we divide by $(z_1+1)$
the rows from the first  to the $(n_1-1)$-th, and by $(z_2+1)$ the remaining ones. We then get
\begin{equation*}
D_*(z_1,z_2)=(z_1+1)^{n_1-1}(z_2+1)^{n_2-1}\det\left[\begin{array}{c|c}
I_{n_1-1}-\zeta_1\mathcal A_{11}&-{{\zeta_1}}\mathcal A_{12}\\
[1ex]\hline\\[-2ex]
-\zeta_2\mathcal A_{12}^{\top}&I_{n_2-1}-\zeta_2\mathcal A_{22}\\[1ex]
\end{array}\right].
\end{equation*}
Comparing the two expressions of $D_*$ yields
\begin{equation}\label{Dred}
\det\left[\begin{array}{c|c}
I_{n_1-1}-\zeta_1\mathcal A_{11}&-{{\zeta_1}}\mathcal A_{12}\\
[1ex]\hline\\[-2ex]
-\zeta_2\mathcal A_{12}^{\top}&I_{n_2-1}-\zeta_2\mathcal A_{22}\\[1ex]
\end{array}\right]=(1-\tau_1\zeta_1^2)^{k_1}(1-\tau_2\zeta_2^2)^{k_2}.
\end{equation}

A direct computation with the left-hand side of~\eqref{Dred}, together
with the fact that $\text{tr}\mathcal A_{\ell\ell}=0$, gives that the
coefficient of $\zeta_1\zeta_2$ is
$$
\sum_{i=1}^{n_1-1}\sum_{j=1}^{n_2-1}\det \left[\begin{array}{cc}
(\mathcal A_{11})_{ii}&(\mathcal A_{12})_{ij}\\
(\mathcal A_{12})_{ij}&(\mathcal A_{22})_{jj}
\end{array}\right]=\text{tr}\mathcal A_{11}\text{tr}\mathcal A_{22}-|\mathcal A_{12}|^2=-|\mathcal A_{12}|^2.
$$

On the other hand, the right-hand side says that such coefficient must vanish.
We conclude that $$\mathcal A_{12}=0,$$ as we wished to prove.

\medskip
\noindent
(ii).
Now we are assuming
$\epsilon_1=0\neq\epsilon_2$. Starting from \eqref{A=0}, let
$$
\tilde B(r,z_2):=\left[\begin{array}{c|c}
I_{n_1}-r\tilde{\mathcal A}_{00}&-r\mathcal A_{02}\\[1ex]
\hline\\[-2ex]
-\frac{z_2-1}{\sqrt{\tau_2}}\mathcal A_{02}^{\top}&(z_2+1)I_{n_2-1}-\frac{z_2-1}{\sqrt{\tau_2}}\mathcal A_{22}
\end{array}\right]
$$
and
$$
\tilde D(r,z_2):=\det\tilde B(r,z_2).
$$
The analogues  of \eqref{DvsD}--\eqref{eqD3} lead to
\begin{eqnarray*}
\tilde D(r,z_2)&=&2^{n_2-1}z_2^{k_2}\det(I_{n_1}-r\tilde{\mathcal A}_{00})\\
&=&2^{n_2-1}(z_2+1)^{n_2-1}\det(I_{n_1}-r\tilde{\mathcal A}_{00})
\det\left(I_{n_2-1}-\frac{z_2-1}{\sqrt{\tau_2}(z_2+1)}\tilde{\mathcal A}_{22}\right),
\end{eqnarray*}
where, this time,
$$
\tilde{\mathcal A}_{22}=\mathcal A_{22}+r\mathcal A_{02}^{\top}(I_{n_1}-r\tilde{\mathcal A}_{00})^{-1}\mathcal A_{02}.
$$

Arguing as the previous case, we can deduce that the eigenvalues of
$\tilde{\mathcal A}_{22}$ are $\pm\sqrt{\tau_2}$, each one with
multiplicity $k_2$. Moreover, for $r$ close to $0$, we have
$$
0=\text{tr}\tilde{\mathcal A}_{22}=\text{tr}\mathcal A_{22}+r|\mathcal A_{02}|^2+O(r^2),
$$
which implies that $$\text{tr}\mathcal A_{22}=0, \qquad \mathcal A_{02}=0.$$

From this point on, we proceed exactly as in the previous case
and reach the desired conclusion.
\end{proof}

\begin{remark}\label{rmkDpol}
On $\mathcal R$, Lemma~\ref{lem:parity} and Lemma~\ref{D-pol} both  fail.
Indeed, for the minimal bi-horospherical examples (see Example~\ref{minimal-bi-horo}),
neither $n_1$ nor $n_2$  is necessarily odd. Besides, for these hypersurfaces,
$\text{tr}\mathcal A_{11}\neq 0\neq \text{tr}\mathcal A_{22}$.
\end{remark}

The next proposition shows that none of the matrices found in Lemma~\ref{D-pol} can be the
shape operator of a hypersurface of \,$\qq$.

\begin{proposition}\label{NoPol}
Assume  $\epsilon_1\le 0$ and $\epsilon_2<0$. Then there is
no hypersurface $\widehat\Sigma\subset\qq$ with angle function $\Theta^2<1$,  whose shape operator blocks
$\mathcal A_{k\ell}$, $0\le k,\ell\le 2$,  satisfy:
\begin{equation}  \label{eq:conditions}
\left\{
\begin{array}{llll}
\mathcal A_{\ell\ell}^2=\tau_\ell I_{n_\ell-1}, &
\mathcal A_{01}=\mathcal A_{02}=\mathcal A_{12}=0, &  {\rm tr}\mathcal A_{\ell\ell}=0, & \text{if} \,\,\, \epsilon_1<0,\\[1ex]
\mathcal A_{22}^2=\tau_2 I_{n_2-1}, & \mathcal A_{02}=0, &  {\rm tr}\,\mathcal A_{22}=0, &  \text{if} \,\,\, \epsilon_1=0.
\end{array}
\right.
\end{equation}
\end{proposition}
\begin{proof}
Assume, by contradiction, that such a hypersurface
$$\widehat\Sigma\subset\qq$$
exists. Since $\Theta^2<1$, setting $N=N_1+N_2$ for
the unit normal to $\widehat\Sigma$,  we have $N_1, N_2\neq0$. Hence the distribution
$$\mathfrak D_\ell:=\mathfrak X(\widehat\Sigma)\cap\mathfrak X(\mathbb Q^{n_{\ell}}_{\epsilon_{\ell}}), \qquad \ell\in\{1,2\},$$
has constant rank $n_{\ell}-1\geq1$. In addition,
being the intersection of $\mathfrak X(\widehat\Sigma)$ with the integrable
distribution of the $\ell$-th factor, $\mathfrak D_{\ell}$ is integrable.

Under the conditions~\eqref{eq:conditions}, the shape operator $A$ of
$\widehat\Sigma$ becomes
\[
A=\left[\begin{array}{c|c|c}
\mathcal A_{00}&0&0\\[1ex]
\hline\\[-2ex]
0&\mathcal A_{11}&0\\
[1ex]\hline\\[-2ex]
0&0&\mathcal A_{22}\\[1ex]
\end{array}\right] \,\, (\epsilon_1,\epsilon_2<0), \qquad
A=\left[\begin{array}{c|c}
\widetilde{\mathcal A}_{00}&0\\[1ex]
\hline\\[-2ex]
 0&\mathcal A_{22}\\[1ex]
\end{array}\right] \,\, (\epsilon_1=0>\epsilon_2).
\]

In any case, the distribution $\mathfrak D_2$ is $A$-invariant, i.e., $A(\mathfrak D_2)\subset\mathfrak D_2$,
and  $AU_1$ is orthogonal to $\mathfrak D_2$. Since $U_1=T/\|T\|$ and $A$ is symmetric, this last condition implies that
$AT$ is orthogonal to $\mathfrak D_2,$
i.e.,
\begin{equation} \label{eq:ATX}
\langle AT,X\rangle=0\qquad \forall X\in\mathfrak D_2.
\end{equation}

Let $\mathcal L\subset \h^{n_{2}}_{\epsilon_{2}}$ be a leaf of
$\mathfrak D_2$ with unit normal
\[
\xi_2=\frac{N_2}{\|N_2\|}\in\mathfrak X(\h^{n_{2}}_{\epsilon_{2}}).
\]

With this orientation, the shape operator $S$
of $\mathcal L$, as a hypersurface of $\h_{\epsilon_2}^{n_2}$, is
\begin{equation} \label{eq:S}
S=\frac{1}{\|N_2\|}A|_{\mathfrak D_2}.
\end{equation}
Indeed, given $X\in\mathfrak X(\mathcal L)$,
using $\nabla\theta=-2AT$,~\eqref{eq:ATX} yields
\[
X(\theta)=\langle\nabla\theta,X\rangle=\langle -2AT,X\rangle=0,
\]
which implies that $\theta$, and so $\|N_2\|^2=(1-\theta)/2$,
is constant on $\mathcal L$. Thus, setting $\xi_1=N_1/\|N_1\|$,
we have $N=\|N_1\|\xi_1+\|N_2\|\xi_2$, and so
\[
AX=-\overbar\nabla_X N=-X(\|N_1\|)\xi_1-\|N_1\|\overbar\nabla_X\xi_1-\|N_2\|\overbar\nabla_X\xi_2.
\]
 But since $\mathfrak D_2$ is $A$-invariant, we must have
 $-X(\|N_1\|)\xi_1-\|N_1\|\overbar\nabla_X\xi_1=0$
 and
 $$AX=-\|N_2\|\overbar\nabla_X\xi_2=\|N_2\|SX,$$
 which gives~\eqref{eq:S}.

 From the conditions~\eqref{eq:conditions}, we have
 $\mathcal A_{22}^2=\tau_2I_{n_2-1}$ and ${\rm tr}\, \mathcal A_{22}=0$. This, together
 with~\eqref{eq:S}, implies that  the eigenvalues of
 $S$  are
 $$\lambda_{\pm}:=\pm\frac{\sqrt{\tau_2}}{\|N_2\|}=\pm\sqrt{-\epsilon_2},$$ each with
multiplicity $(n_2-1)/{2}$. In particular, $n_2$ is odd.

Therefore  $\mathcal L$ is a hypersurface of the
hyperbolic space $\h_{\epsilon_2}^{n_2}$
with  constant principal curvatures. In this context, Cartan's formula
applies (cf. \cite{cartan1,cartan3}):
\[
\sum_{j\neq i}m_j\frac{\epsilon_2+\lambda_i\lambda_j}{\lambda_i-\lambda_j}=0,
\]
where, for each index $j$,  $\lambda_j$ is  a principal curvature  of the hypersurface
with multiplicity $m_j$.
Specializing to our case, this formula reduces to
$$0=\epsilon_2+\lambda_+\lambda_-=2\epsilon_2<0,$$
giving a contradiction.
\end{proof}


{\subsection{The proof.}
Now, we are in a position to prove Theorem~\ref{thm-isoparametric-->thetaconstant}, which we
restate for the reader's convenience.

\isoimpliesthetact*
\begin{proof}
After rearranging the factors and restricting to a suitable open set, we can assume
 $\epsilon_2\neq 0$ and $\Theta^2<1$ on $\Sigma.$ We separate the proof into two cases.

\medskip
\noindent \emph{Case 1:} ${D(r)}$ \emph{is not a polynomial.} By Lemma~\ref{lem:polynomial-equivalence},
we have $\Lambda\neq\{0\}$. Let
$$\lambda^*\in\Lambda\setminus\{0\}.$$

Given $p\in\Sigma,$ because $\Lambda$ is independent of $p$,  there exist
$a_\ell=a_\ell(p)\in\mathcal N_\ell$, $\ell\in\{1,2\}$,  such that
\begin{equation}\label{l*}
\lambda^*=a_1(p)\sqrt{\tau_1(p)}+a_2(p)\sqrt{\tau_2(p)}, \qquad (a_1(p),a_2(p))\ne(0,0),
\end{equation}
where the second condition comes from $\lambda^*\ne0$.

Setting $\nu=\|N_1\|\in(0,1)$, we have
$\|N_2\|=\sqrt{1-\nu^2}>0$. Then, considering that
$\tau_\ell=-\epsilon_\ell\|N_\ell\|^2$, it follows from~\eqref{l*}
that, at every point of $\Sigma$, the value $\nu(p)$ satisfies one of the finitely
many equations
\begin{equation}\label{eq:alpha}
\alpha_1t+\alpha_2\sqrt{1-t^2}=\lambda^*,
\qquad
\alpha_\ell=a_\ell\sqrt{-\epsilon_\ell},
\qquad
(a_1,a_2)\in\mathcal N_1\times\mathcal N_2\setminus\{(0,0)\}.
\end{equation}

For each fixed pair $(a_1,a_2)$, equation~\eqref{eq:alpha} implies
\[
(\alpha_1^2+\alpha_2^2)t^2-2\alpha_1\lambda^*t+(\lambda^*)^2-\alpha_2^2=0.
\]

Hence, if $\alpha_1^2+\alpha_2^2\ne0$,
$\nu(p)$ is a solution of the quadratic equation~\eqref{eq:alpha}.
If $\alpha_1^2+\alpha_2^2=0$, then both
$\alpha_1$ and $\alpha_2$ are nonzero, for
$(\alpha_1,\alpha_2)\ne(0,0)$. Since $\lambda^*\ne 0$,
it implies that $\nu(p)$
is a solution of the  affine equation:
\[
-2\alpha_1\lambda^*t+(\lambda^*)^2-\alpha_2^2=0.
\]

Since $\mathcal N_1\times\mathcal N_2$ is finite, we
conclude that $\nu(\Sigma)$ is contained in a finite set, and so $\nu$ is constant, for
it is continuous and $\Sigma$ is connected.
Consequently, the same is true for
\[
\Theta=2\|N_1\|^2-1=2\nu^2-1.
\]

\noindent
\emph{Case 2}: ${D(r)}$ \emph{is a polynomial.}
Suppose that $\widehat\Sigma:=\Sigma\setminus\mathcal R\ne\emptyset.$
Since $\widehat\Sigma$ is open in $\Sigma$ (cf.~Rem.~\ref{rem:Rclosed}), it
is a hypersurface of $\qq$. Moreover, by Lemma~\ref{D-pol}, $\epsilon_1\le0$, $\epsilon_2<0$, and
$\widehat\Sigma$ satisfies
the conditions~\eqref{eq:conditions}. But, by Proposition~\ref{NoPol}, this is impossible.

Therefore $\Sigma=\mathcal R,$ that is, $\big(\sqrt{\tau_1(p)},\sqrt{\tau_2(p)}\big)$
is $(n_1,n_2)$-dependent for all $p\in\Sigma$. By Lemma~\ref{lemma-n1n2dep}(iii), for each $\ell\in\{1,2\},$
there exists $a_\ell(p)\in\{1,\dots, n_\ell-1\}$ such that
\begin{equation} \label{eq:a1}
 a_1^2(p){\tau_1(p)}=a_2^2(p){\tau_2(p)}, \qquad \big(a_1(p),a_2(p)\big)\ne (0,0),
\end{equation}
where the second condition comes from the $(n_1,n_2)$-dependence of $\big(\sqrt{\tau_1(p)},\sqrt{\tau_2(p)}\big)$.
We also have  $\tau_1\tau_2>0$ (cf.~Rem.~\ref{rem:Rclosed}).
Thus $a_1(p)$ and $a_2(p)$ are both nonzero, and it
follows from~\eqref{eq:a1} that the continuous function $\tau_1/\tau_2$ takes values in the
finite set
\[
\left\{\frac{a_2^2}{a_1^2}\,;\,
1\leq a_1\leq n_1-1,\quad1\leq a_2\leq n_2-1\right\}.
\]

Therefore $\tau_1/\tau_2$ is constant, for $\Sigma$ is connected. Moreover,  $\tau_\ell=-\epsilon_\ell\|N_\ell\|^2$
gives
\[
\frac{\tau_1}{\tau_2}
=\frac{\epsilon_1}{\epsilon_2}\frac{\nu^2}{1-\nu^2}, \qquad \nu^2=\|N_1\|^2=\frac{\theta+1}{2}\in(0,1),
\]
which implies that $\nu^2$, and then $\theta$, is constant on $\Sigma$. This completes our proof.
\end{proof}

\begin{remark} \label{rem:gaps01}
Theorem 1.1 of~\cite{tan02} states that every isoparametric hypersurface in
$\h^n\times\R^m$ or $\s^n\times\R^m$ has constant angle function. In its proof, the
authors employ a technique that we introduced in~\cite{dLP} and that, as in the present
paper, involves the function $D(r)$. However, several of their arguments are not justified,
and some are, in fact, incorrect. For instance, they conclude from their equation~(4.5)
that $D(r)$ is a polynomial near $r=0$ and then claim that this is impossible because $D(r)$
cannot be a polynomial. As we have seen here, this claim is false: the polynomiality of
$D(r)$ is not algebraically impossible; rather, it must be ruled out by geometric arguments.
\end{remark}


\section{Proof of Theorem~\ref{thm-classification01}}  \label{sec-classification01}

Let $\Sigma\subset\qq$ be a connected isoparametric hypersurface
with constant angle function $\Theta\in(-1,1)$.
Recall that, in this case, we have  $T\neq0$. In addition,
\[
 \|N_1\|^2=\frac{1+\Theta}{2},
 \qquad
 \|N_2\|^2=\frac{1-\Theta}{2},
 \qquad\text{and}\qquad
 \tau_\ell=-\epsilon_\ell\|N_\ell\|^2, \,\,\, \ell\in\{1,2\},
\]
are constants, with both $\|N_1\|$ and $\|N_2\|$ positive.

We will use the frame
\[
\{U^1,\ldots,U^{n-1}\}\subset\mathfrak X(\Sigma), \qquad U_1=\frac{T}{\|T\|},
\]
considered in Section~\ref{sec-constantangle}, and
the corresponding blocks $\mathcal A_1,\mathcal A_2,
\mathcal A_3$ of the shape operator defined there.

Finally, recall that
\[
 c_\tau(r)=
 \begin{cases}
  \cosh(\sqrt\tau r),&\tau>0,\\
  1,&\tau=0,\\
  \cos(\sqrt{-\tau}r),&\tau<0,
 \end{cases}
 \qquad\quad
 s_\tau(r)=
 \begin{cases}
  \dfrac{\sinh(\sqrt\tau r)}{\sqrt\tau},&\tau>0,\\[.5ex]
  r,&\tau=0,\\[1ex]
  \dfrac{\sin(\sqrt{-\tau}r)}{\sqrt{-\tau}},&\tau<0,
 \end{cases}
\]
and set, for  $c_\tau(r)\ne 0,$
\[
\mathcal Z_\tau(r)=\frac{s_\tau(r)}{c_\tau(r)}\cdot
\]

In this setting, assuming further that
\[
\epsilon_1\epsilon_2\leq0,
\]
we have the following $4$ lemmas.

\begin{lemma}\label{lem:independencez1z2}
The functions $\mathcal Z_{\tau_1}$ and $\mathcal Z_{\tau_2}$ are algebraically
independent over $\R$.
\end{lemma}

\begin{proof}
Suppose, to the contrary, that there is a nonzero polynomial
\[
 \mathcal P(x,y)=\sum_{j=0}^{d}p_j(x)y^j,
 \qquad p_d\neq0,
\]
such that
\begin{equation}\label{eq:algebraic-relation}
 \mathcal P\bigl(\mathcal Z_{\tau_1}(r),\mathcal Z_{\tau_2}(r)\bigr)=0
\end{equation}
on an open interval containing zero.

After interchanging the indices, we have the following three cases:
\[
 \tau_2<0<\tau_1,
 \qquad
 \tau_2<0=\tau_1,
 \qquad
 \tau_1=0<\tau_2.
\]

\medskip
\noindent
\emph{ Case 1:} ${\tau_2<0<\tau_1}$.
Set $a=\sqrt{\tau_1}$ and $b=\sqrt{-\tau_2}$, so that
\[
 \mathcal Z_{\tau_1}(r)=\frac{\tanh(ar)}a
 \qquad\textup{and}\qquad
 \mathcal Z_{\tau_2}(r)=\frac{\tan(br)}b\cdot
\]
Consider the function
\begin{align*}
 G(r)
 &:=\cos^d(br)
 \mathcal P\left(\frac{\tanh(ar)}a,\frac{\tan(br)}b\right)\\
 &=\sum_{j=0}^{d}\frac1{b^j}
 p_j\left(\frac{\tanh(ar)}a\right)
 \sin^j(br)\cos^{d-j}(br).
\end{align*}
The second expression shows that $G$ is real analytic on the whole real
line: multiplication by $\cos^d(br)$ has removed the denominators that
occur in the powers of $\tan(br)$. Since $G$ vanishes on an open
interval, it vanishes identically. In particular, at the points
\[
 \rho_k=\frac{\pi/2+k\pi}{b},
 \qquad k\in\mathbb Z,
\]
we have
\[
 0=G(\rho_k)
 =\frac{\sin^d(b\rho_k)}{b^d}
 p_d\left(\frac{\tanh(a\rho_k)}a\right),
\]
and hence
\[
 p_d\left(\frac{\tanh(a\rho_k)}a\right)=0
 \qquad \forall k\in\mathbb Z.
\]

So, for $k\geq0$, the numbers
${\tanh(a\rho_k)}/a$
are pairwise distinct, because the function
$r\mapsto\tanh(ar)/a$ is strictly increasing. Thus the nonzero
polynomial $p_d$ has infinitely many distinct zeros; a contradiction.

\medskip
\noindent
\emph{Case 2:} ${\tau_2<0=\tau_1}$. Now $\mathcal Z_{\tau_1}(r)=r$. Setting $b=\sqrt{-\tau_2}$, the preceding
argument, with $\tanh(ar)/a$ replaced by $r$, shows that the
real analytic function
\[
 G(r)=\cos^d(br)\mathcal P\left(r,\frac{\tan(br)}b\right)
\]
vanishes identically. Evaluating at
$\rho_k=(\pi/2+k\pi)/b$ gives
\[
 p_d(\rho_k)=0
 \qquad \forall k\in\mathbb Z,
\]
again contradicting $p_d\neq0$.

\medskip
\noindent
\emph{Case 3:} ${\tau_1=0<\tau_2}$.
Set $a=\sqrt{\tau_2}$. The left-hand side of
\[
 \mathcal P\left(r,\frac{\tanh(ar)}a\right)=0
\]
is real analytic on $\mathbb R$. Hence the identity, initially known on
an open interval, holds for every $r\in\mathbb R$. Put $u=e^{-2ar}$.
Since
\[
 \frac{\tanh(ar)}a=\frac{1-u}{a(1+u)}\,,
\]
multiplication by $(1+u)^d$ gives
\begin{equation}\label{eq:exponential-polynomial}
 \sum_{k=0}^{d}q_k(r)u^k=0,
\end{equation}
where $q_0,\ldots,q_d\in\mathbb R[r]$ are not all zero. Indeed, the
change of variables
\[
 y=\frac{1-u}{a(1+u)}
\]
is invertible, with inverse $u=(1-ay)/(1+ay)$, and therefore cannot
transform a nonzero polynomial into the zero polynomial after clearing
denominators.

Let $k_0$ be the smallest index for which $q_{k_0}\neq0$. Since
$u=e^{-2ar}$, dividing \eqref{eq:exponential-polynomial} by $u^{k_0}$
yields
\[
 q_{k_0}(r)+\sum_{k=k_0+1}^{d}
 q_k(r)e^{-2a(k-k_0)r}=0.
\]
As $r\to+\infty$, every term in the sum tends to zero. Hence
$q_{k_0}(r)\to0$, which is impossible because $q_{k_0}$ is a nonzero
polynomial. This final contradiction proves the lemma.
\end{proof}

\begin{lemma}\label{lem:constant-characteristic-polynomials}
The characteristic polynomials of $\mathcal A_1$ and $\mathcal A_3$
are constant on $\Sigma$.
\end{lemma}

\begin{proof}
For each $p\in\Sigma$, define the polynomial $\Phi_p\in\mathbb R[x,y]$
by
\begin{equation}\label{eq:Phi}
 \Phi_p(x,y)=\det
 \begin{bmatrix}
  I_{n_1-1}-x\mathcal A_1(p)&-x\mathcal A_2(p)\\[1ex]
  -y\mathcal A_2(p)^\top&I_{n_2-1}-y\mathcal A_3(p)
 \end{bmatrix}.
\end{equation}

A direct computation gives that the function
$D(r)$ defined in~\eqref{eq-D&cs} satisfies
\begin{align}
 D(r) &=c_{\tau_1}(r)^{n_1-1}c_{\tau_2}(r)^{n_2-1}
 \Phi_p\bigl(\mathcal Z_{\tau_1}(r),\mathcal Z_{\tau_2}(r)\bigr),
 \label{eq:normal-jacobian}
\end{align}
which implies that $\Phi_p\bigl(\mathcal Z_{\tau_1}(r),\mathcal Z_{\tau_2}(r)\bigr)$
is independent of $p$. Thus, for all
sufficiently small $r$,
\begin{equation}\label{eq:Phi-on-curve}
 \Phi_p\bigl(\mathcal Z_{\tau_1}(r),\mathcal Z_{\tau_2}(r)\bigr)
 =\Phi_q\bigl(\mathcal Z_{\tau_1}(r),\mathcal Z_{\tau_2}(r)\bigr)
 \qquad \forall p,q\in\Sigma.
\end{equation}

Fix $p,q\in\Sigma$ and set
\[
 \mathcal P(x,y)=\Phi_p(x,y)-\Phi_q(x,y).
\]

The equality~\eqref{eq:Phi-on-curve} and the algebraic independence proved in
Lemma~\ref{lem:independencez1z2}  yield $\mathcal P\equiv0$. Consequently,
\[
 \Phi_p(x,y)=\Phi_q(x,y)
 \qquad \forall p,q\in\Sigma.
\]
Setting $y=0$ in this equality gives
\[
 \det\bigl(I_{n_1-1}-x\mathcal A_1(p)\bigr)
 =\det\bigl(I_{n_1-1}-x\mathcal A_1(q)\bigr),
\]
whereas setting $x=0$ gives
\[
 \det\bigl(I_{n_2-1}-y\mathcal A_3(p)\bigr)
 =\det\bigl(I_{n_2-1}-y\mathcal A_3(q)\bigr).
\]
Thus the characteristic polynomials of $\mathcal A_1$ and
$\mathcal A_3$ are constant on $\Sigma$.
\end{proof}


\begin{lemma} \label{lem:spectral-identities}
Let $$\mathcal A_1=\sum_{\lambda\in\spec\mathcal A_1}\lambda P_\lambda
\qquad\text{and}\qquad
\mathcal A_3=\sum_{\mu\in\spec\mathcal A_3}\mu Q_\mu$$
be the spectral decompositions of $\mathcal A_1$ and $\mathcal A_3$, respectively.
Then
\begin{equation} \label{eq:spectral-identities}
  c\|P_\lambda\mathcal A_2\|^2
    =m_\lambda(\lambda^2-\tau_1)
  \qquad\text{and}\qquad
  \|\mathcal A_2Q_\mu\|^2
    =c\,m_\mu(\mu^2-\tau_2),
\end{equation}
where $m_\lambda$ and $m_\mu$ denote the multiplicities of
the eigenvalues $\lambda$ and $\mu$, respectively.
In particular, the following hold:
\begin{itemize}[parsep=.5ex]
\item[\rm (i)] $\lambda^2\geq\tau_1$ for every $\lambda\in \spec \mathcal A_1$, and $\mu^2\geq\tau_2$ for every $\mu\in \spec \mathcal A_3$,
\item[\rm (ii)] if either $\tau_1$ or $\tau_2$ is negative, then $\mathcal A_2\neq0$ at every point of $\Sigma$.
\end{itemize}
\end{lemma}

\begin{proof}
By Lemma~\ref{lem:constant-characteristic-polynomials},
the characteristic polynomials of $\mathcal A_1$ and of $\mathcal A_3$
are constant on $\Sigma$. Hence
$\tr \mathcal A_1^{k+1}$ is constant for every integer $k\ge 0$, and so
$$0=T(\tr \mathcal A_1^{k+1})=(k+1)\tr\bigl(\mathcal A_1^{k}T(\mathcal A_1)\bigr),$$
which implies that
$\tr\bigl(\mathcal A_1^k\widetilde T(\mathcal A_1)\bigr)=0$  {for every integer} $k\geq0.$
Multiplying the first equation of the system~\eqref{eq:system} by $\mathcal A_1^k$ and taking traces then  yields
\begin{equation} \label{eq:spec}
  \sum_{\lambda\,\in \,\spec{\mathcal A_1}}\lambda^k
  \left[m_\lambda(\lambda^2-\tau_1)
    -c\,\tr\bigl(P_\lambda \mathcal A_2A_2^{\top}\bigr)\right]=0
  \qquad \forall k\geq0.
\end{equation}

Let $\spec \mathcal A_1=\{\lambda_1,\ldots,\lambda_r\}$
and set
\[
b_\alpha=m_{\lambda_\alpha}(\lambda_\alpha^2-\tau_1)-c\operatorname{tr}\!\left(
          P_{\lambda_\alpha}\mathcal A_2\mathcal A_2^\top\right),\qquad 1\leq\alpha\leq r.
\]

For each $\lambda=\lambda_\ell\in\spec{\mathcal A_1},$ take
$j=0,1,\ldots,r-1$ in~\eqref{eq:spec}. This gives
\[
Vv^\top=0, \qquad v=(b_1,\dots, b_r),
\]
where $V$ is the standard Vandermonde matrix of order $r$ determined by
the pairwise distinct eigenvalues $\lambda_1,\dots, \lambda_r$.
In particular,
\[
\det V=\prod_{1\leq\alpha<\beta\leq r}
       (\lambda_\beta-\lambda_\alpha)\ne0,
\]
which implies that $v=(b_1,\dots, b_r)=0$. Therefore
\[
c\operatorname{tr}\!\left(P_\lambda\mathcal A_2\mathcal A_2^\top\right)=m_\lambda(\lambda^2-\tau_1)
  \qquad \forall\lambda\in\spec\mathcal A_1.
\]

Now, since $P_\lambda$ is an orthogonal projection, we have $P_\lambda^2=P_\lambda,$ so
\[
\operatorname{tr}\!\left(
  P_\lambda\mathcal A_2\mathcal A_2^\top
\right)
 =\operatorname{tr}\!\left(
   \mathcal A_2^\top P_\lambda\mathcal A_2
 \right) \\
 =\|P_\lambda\mathcal A_2\|^2,
\]
which proves the first equality in~\eqref{eq:spectral-identities}.

An analogous argument using the second equation of~\eqref{eq:system}
gives the second equality in~\eqref{eq:spectral-identities}.
Finally, the assertions (i) and (ii) follow
directly from~\eqref{eq:spectral-identities}.
\end{proof}


The classification of isoparametric hypersurfaces of $\qq$ when
$\epsilon_1\epsilon_2\le 0$ will be a consequence of the  next
two propositions. The first one gives conditions under which an
isoparametric hypersurface of $\qq$ has constant angle equal to $\pm1$,
whereas the second shows that an isoparametric hypersurface of
$\R^{n_1}\times\h_{\epsilon_2}^{n_2}$ whose principal directions split properly
at every point is necessarily an open set of a flat-horospherical hypersurface.

\begin{proposition}\label{prop:theta-one}
Let
$\Sigma\subset\mathbb Q_{\epsilon_1}^{n_1}\times\mathbb Q_{\epsilon_2}^{n_2}$ be a
connected isoparametric hypersurface, where
$\epsilon_1\epsilon_2\leq0$. Then $\Theta^2\equiv1$ on $\Sigma$,
provided that one of the following occurs:
\begin{itemize}[parsep=.7ex]
 \item[\rm (i)] $\epsilon_i>0$ for some $i\in\{1,2\}$;
 \item[\rm (ii)] $\epsilon_1\epsilon_2=0$ with some $\epsilon_i<0$, and
 $\mathcal A_2\neq0$ at some $p\in\Sigma$.
\end{itemize}
\end{proposition}

\begin{proof}
By Theorem~\ref{thm-isoparametric-->thetaconstant},
$\Sigma$  has constant angle $\Theta$. We will
assume $\Theta^2<1$ and derive a contradiction. Up to a change of
factors, we have three cases:
\[
 \h_{\epsilon_1}^{n_1}\times\mathbb S_{\epsilon_2}^{n_2},
 \qquad
 \R^{n_1}\times\mathbb S_{\epsilon_2}^{n_2},
 \qquad
 \h_{\epsilon_1}^{n_1}\times\R^{n_2},
\]
which we shall consider separately.

\medskip
\noindent
\emph {Case 1:} $\h_{\epsilon_1}^{n_1}\times\mathbb S_{\epsilon_2}^{n_2}$.
In this case,
\[
 \tau_2<0<\tau_1,
\]
and so, by Lemma~\ref{lem:spectral-identities}, $\mathcal A_2$ never
vanishes on $\Sigma$.

As we have seen in the proof of Proposition~\ref{NoPol}, the distribution
\[
\mathfrak D_1=\mathfrak X(\Sigma)\cap\mathfrak X(\mathbb H_{\epsilon_1}^{n_1})
\]
is integrable, and  a leaf
$\mathcal L_p\subset\h_{\epsilon_1}^{n_1}$ of $\mathfrak D_1$ has shape operator
$\mathcal A_1/\|N_1\|$. Hence,  by
Lemma~\ref{lem:constant-characteristic-polynomials}, $\mathcal L_p$ has
constant principal curvatures and is therefore an isoparametric hypersurface of
$\h_{\epsilon_1}^{n_1}$.

Writing $\delta=\sqrt{|\epsilon_1|}$ and
$\kappa=\lambda/\|N_1\|$, Lemma~\ref{lem:spectral-identities} gives
\[
 \kappa^2=\frac{\lambda^2}{\|N_1\|^2}
 \geq\frac{\tau_1}{\|N_1\|^2}
 =|\epsilon_1|=\delta^2,
\]
so that $|\kappa|\geq\delta>0$ for every principal curvature $\kappa$
of $\mathcal L_p$.

Recall that the isoparametric hypersurfaces of
$\h_{\epsilon_1}^{n_1}$ are: totally geodesic hypersurfaces
($\kappa=0$), equidistant hypersurfaces
($\kappa=\delta\tanh(\delta d)$, $|\kappa|<\delta$), horospheres
($|\kappa|=\delta$), geodesic spheres
($\kappa=\delta\coth(\delta d)$, $|\kappa|>\delta$), and tubes about a
totally geodesic $\h_{\epsilon_1}^j$, $j<n_1-1$, whose principal
curvatures are $\delta\coth(\delta d)$ and
$\delta\tanh(\delta d)$; see~\cite[pp. 97--98]{ryan}.
The condition $|\kappa|\geq\delta$ rules out
all the non-umbilical ones, and also the totally geodesic and equidistant
ones. Hence $\mathcal L_p$ is umbilical:
\[
 \mathcal A_1=\lambda I_{n_1-1},
 \qquad
 \lambda^2\geq\tau_1,
 \qquad
 \lambda=\text{constant}.
\]
In particular, $\widetilde T(\mathcal A_1)=0$, and the first equation of
\eqref{eq:system} becomes
\[
 \mathcal A_2\mathcal A_2^\top=\sigma^2I_{n_1-1},
 \qquad
 \sigma^2:=\frac{\lambda^2-\tau_1}{c}\geq0.
\]
Since $\mathcal A_2\neq0$, we have $\sigma^2>0$. Set
\[
 W=\frac1\sigma\mathcal A_2, \qquad\text{so that}
 \,\,
 WW^\top=I_{n_1-1}.
\]

We now show that
\begin{equation}\label{eq:WA3W00}
 W\mathcal A_3W^\top=\frac\lambda cI_{n_1-1}.
\end{equation}
Since $\mathcal A_2\mathcal A_2^\top$ is constant,
\begin{equation}\label{eq:derivative-A2A2}
 \widetilde T(\mathcal A_2)\mathcal A_2^\top
 +\mathcal A_2\widetilde T(\mathcal A_2)^\top=0.
\end{equation}
The third equation of \eqref{eq:system} yields
\begin{equation}\label{eq:derivative-A2}
 \widetilde T(\mathcal A_2)
 =\lambda\mathcal A_2-c\mathcal A_2\mathcal A_3.
\end{equation}
Substituting \eqref{eq:derivative-A2} into
\eqref{eq:derivative-A2A2} and using $\mathcal A_2=\sigma W$ gives
\eqref{eq:WA3W00}.

Next, we show that
\begin{equation}\label{eq:WA3sqW}
 W\mathcal A_3^2W^\top
 =\frac1{3c}\left(\frac{2\lambda^2}{c}
                   +\sigma^2+c\tau_2\right)I_{n_1-1}.
\end{equation}
By~\eqref{eq:WA3W00}, the endomorphism
$\mathcal A_2\mathcal A_3\mathcal A_2^\top$ is constant. Hence
\begin{equation}\label{eq:derivative-A2A3A2}
 \widetilde T(\mathcal A_2)\mathcal A_3\mathcal A_2^\top
 +\mathcal A_2\widetilde T(\mathcal A_3)\mathcal A_2^\top
 +\mathcal A_2\mathcal A_3\widetilde T(\mathcal A_2)^\top=0.
\end{equation}
The second equation of \eqref{eq:system} gives
\begin{equation}\label{eq:derivative-A3}
 \widetilde T(\mathcal A_3)
 =\mathcal A_2^\top\mathcal A_2
  -c\bigl(\mathcal A_3^2-\tau_2I_{n_2-1}\bigr).
\end{equation}
Applying \eqref{eq:derivative-A2} and \eqref{eq:derivative-A3} to
\eqref{eq:derivative-A2A3A2} yields
\[
 \frac{2\lambda^2\sigma^2}{c}I_{n_1-1}
 -3c\mathcal A_2\mathcal A_3^2\mathcal A_2^\top
 +\sigma^4I_{n_1-1}
 +c\tau_2\sigma^2I_{n_1-1}=0,
\]
which implies \eqref{eq:WA3sqW}.

Finally, set
\[
 \varPi=W^\top W,
 \qquad
 Z=(I_{n_2-1}-\varPi)\mathcal A_3W^\top.
\]
Since $WW^\top=I_{n_1-1}$, one has $\varPi^2=\varPi$; hence $\varPi$ is an
orthogonal projection on $\mathbb R^{n_2-1}$. Therefore
\[
 \begin{split}
 W\mathcal A_3^2W^\top-(W\mathcal A_3W^\top)^2
 &=W\mathcal A_3(I_{n_2-1}-\varPi)\mathcal A_3W^\top\\
 &=Z^\top Z.
 \end{split}
\]
Since $Z^\top Z$ is positive semi-definite,
\eqref{eq:WA3W00} and \eqref{eq:WA3sqW} give
\[
 \frac1{3c}\left(\frac{2\lambda^2}{c}
                  +\sigma^2+c\tau_2\right)
 \geq\frac{\lambda^2}{c^2}\cdot
\]
Substituting $c\sigma^2=\lambda^2-\tau_1$, this reduces to
$c^2\tau_2\geq\tau_1$, which is impossible because
$\tau_2<0<\tau_1$.

\medskip
\noindent
\emph{Case 2:} $\R^{n_1}\times\mathbb S_{\epsilon_2}^{n_2}$.
Now, we have
\[
 \tau_2<0=\tau_1,
\]
and again Lemma~\ref{lem:spectral-identities} gives
$\mathcal A_2\neq0$ on $\Sigma$.

As in the previous case, a leaf $\mathcal L_p\subset\mathbb R^{n_1}$ has
constant principal curvatures and is therefore an open set of a
hyperplane, a sphere, or a tube about an affine subspace of codimension
at least two. The hyperplane case $\mathcal A_1=0$ is
impossible, because the first equation of~\eqref{eq:system} would give
$\mathcal A_2=0$. Consequently, there exists a constant eigenvalue $\lambda\ne 0$ such that
\[
 \mathcal A_1=\lambda P_\lambda,
\]
where $P_\lambda$ is the orthogonal projection onto the  $\lambda$-eigenspace.

Let $P_0$ be the orthogonal projection onto the zero eigenspace. Then we have
\[
P_0=I_{n_1-1}-P_\lambda.
\]

By the first identity
in~\eqref{eq:spectral-identities}, applied to the zero eigenspace,
$\|P_0\mathcal A_2\|^2=0.$
Thus
\begin{equation} \label{eq:A1=lP}
\mathcal A_2=P_\lambda\mathcal A_2.
\end{equation}
Moreover,  differentiating $P_\lambda^2=P_\lambda$ along
$\widetilde T$ gives
\[
\widetilde T(P_\lambda)P_\lambda+P_\lambda\widetilde T(P_\lambda)=\widetilde T(P_\lambda).
\]
Hence, multiplying by $P_\lambda$ on the left, we have
\[
P_\lambda\widetilde T(P_\lambda)P_\lambda=0.
\]

On the other hand, the first equation of~\eqref{eq:system} reads
\begin{equation} \label{eq:lambdaT}
 \lambda\widetilde T(P_\lambda)
 =\lambda^2P_\lambda-c\mathcal A_2\mathcal A_2^\top.
\end{equation}
By~\eqref{eq:A1=lP}, we have $\mathcal A_2\mathcal A_2^\top=P_\lambda\mathcal A_2\mathcal A_2^\top P_\lambda$.
Therefore, multiplying~\eqref{eq:lambdaT} on the left and on the right by $P_\lambda$ gives
\[
\mathcal A_2\mathcal A_2^\top=\sigma^2P_\lambda,
 \qquad
 \sigma^2:=\frac{\lambda^2}{c}>0.
\]
Substituting this equality back into~\eqref{eq:lambdaT} implies that
$\widetilde T(P_\lambda)=0$, and so $\mathcal A_2\mathcal A_2^\top$ is constant
on $\Sigma$.

Set $W=\mathcal A_2/\sigma$. Then
\[
 WW^\top=P_\lambda,
 \qquad
 P_\lambda W=W.
\]
The computations in \eqref{eq:derivative-A2A2}--\eqref{eq:derivative-A3},
with $I_{n_1-1}$ replaced by $P$, now give
\begin{equation} \label{eq:WA3W}
 W\mathcal A_3W^\top=\frac\lambda cP_\lambda,
 \qquad
 W\mathcal A_3^2W^\top
 =\frac1{3c}\left(\frac{2\lambda^2}{c}
                   +\sigma^2+c\tau_2\right)P_\lambda.
\end{equation}

Proceeding then as in the previous case, we conclude that
\[
 W\mathcal A_3^2W^\top-(W\mathcal A_3W^\top)^2
\]
is positive semi-definite. Therefore, using~\eqref{eq:WA3W} and
restricting to the image of $P_\lambda$ yields
$$c^2\tau_2\geq\tau_1=0,$$
contradicting $\tau_2<0$.

\medskip
\noindent
\emph{Case 3:} $\h_{\epsilon_1}^{n_1}\times\R^{n_2}$.
In this case,
\[
 \tau_2=0<\tau_1.
\]
Since the first factor is hyperbolic, the reasoning from Case 1 gives
\[
 \mathcal A_1=\lambda I_{n_1-1},
 \qquad
 \lambda^2\geq\tau_1,
 \qquad
 \lambda=\text{constant}.
\]
The first equation of \eqref{eq:system} then gives
\[
 \mathcal A_2\mathcal A_2^\top=\sigma^2I_{n_1-1},
 \qquad
 \sigma^2:=\frac{\lambda^2-\tau_1}{c}\geq0.
\]
By hypothesis, $\mathcal A_2(p)\neq0$ at some $p\in\Sigma$, so
$\sigma^2>0$. Defining $W=\mathcal A_2/\sigma$ and repeating the
argument of Case 1 gives
\[
 \tau_1\leq c^2\tau_2=0,
\]
which contradicts $\tau_2=0<\tau_1$.

Since the hypothesis $\Theta^2<1$ leads to a contradiction in all three
cases, we conclude that $\Theta^2\equiv1$. This finishes the proof.
\end{proof}


\begin{proposition}\label{prop:flat-horospherical}
Let $\Sigma\subset\R^{n_1}\times\h_{\epsilon_2}^{n_2}$ be a connected
isoparametric hypersurface whose principal directions split properly
everywhere. Then $\Sigma$ is an open set of a flat-horospherical
hypersurface.
\end{proposition}

\begin{proof}
By Theorem~\ref{thm-isoparametric-->thetaconstant}, $\Sigma$ has constant angle.
Since the principal directions split properly at every point of
$\Sigma$, locally, its tangent bundle splits orthogonally as
\[
 \mathfrak X(\Sigma)=\mathfrak D\oplus\mathfrak D_1\oplus\mathfrak D_2,
\]
where $\mathfrak D,\mathfrak D_1$ and $\mathfrak D_2$ are $A$-invariant
distributions on $\Sigma$ such that
\[
 \mathfrak D=\operatorname{span}\{T\},
 \qquad
 \mathfrak D_1\subset\mathfrak X(\R^{n_1}),
 \qquad
 \mathfrak D_2\subset\mathfrak X(\h_{\epsilon_2}^{n_2}).
\]

We first prove that $\mathfrak D_1$ is integrable. Given
$X,Y\in\mathfrak D_1$, the bracket $[X,Y]$ is tangent to the first factor,
and hence
\[
 \langle[X,Y],Z\rangle=0
 \qquad \forall Z\in\mathfrak D_2.
\]
Moreover,
\[
 \langle[X,Y],T\rangle
 =\langle[X,Y],T_1\rangle
 =(1-\Theta)\langle[X,Y],N_1\rangle
 =(1-\Theta)\langle[X,Y],N\rangle=0,
\]
where in the last equality we used $[X,Y]\in\mathfrak X(\Sigma)$.
Therefore $[X,Y]\in\mathfrak D_1$, and $\mathfrak D_1$ is integrable.
The same argument proves that $\mathfrak D_2$ is integrable.

Now, we prove the compatibility of these distributions with the
$T$-flow. Since $P$ is parallel, $\Theta$ is constant, and
$T=PN-\Theta N$, one has
\[
 \overbar\nabla_XT=(\Theta I-P)AX.
\]
Thus, for every $ X\in\mathfrak D_1$ and $Y\in\mathfrak D_2$,
\[
 \overbar\nabla_XT=(\Theta-1)AX\in\mathfrak D_1
 \qquad\text{and}\qquad
 \overbar\nabla_YT=(\Theta+1)AY\in\mathfrak D_2.
\]

Since $X\in\mathfrak D_1$ if and only if $PX=X$, and $P$ is parallel, we have
\[
P\overbar\nabla_TX=\overbar\nabla_TPX=\overbar\nabla_TX \qquad \forall X\in\mathfrak D_1,
\]
that is, $\overbar\nabla_TX\in\mathfrak D_1$ for every $X\in\mathfrak D_1$. Thus
$$[T,X]=\overbar\nabla_TX-\overbar\nabla_XT\in \mathfrak D_1.$$
Similarly,
$[T,Y]\in\mathfrak D_2$ for $Y\in\mathfrak D_2$. Finally, for
$X\in\mathfrak D_1$ and $Y\in\mathfrak D_2$,
\[
 \langle[X,Y],T\rangle
 =-\langle Y,\overbar\nabla_XT\rangle
  +\langle X,\overbar\nabla_YT\rangle=0.
\]
Hence $\mathfrak D_1\oplus\mathfrak D_2$ is integrable, and the $T$-flow
preserves each of the two foliations.

We therefore obtain a local
diffeomorphism
\[
 \Psi:M_1\times M_2\times I\longrightarrow\Sigma
\]
such that the $M_i$-slices are leaves of $\mathfrak D_i$, $i=1,2$, and
the curves
\[
 s\in I\longmapsto\Psi(p_1,p_2,s)\in\Sigma
\]
are integral curves of $T$, after a constant reparameterization.

Denote by $\pi_i$ the projection onto the $i$-th factor and set
\[
 \Phi_i(p_1,p_2,s)=\pi_i\bigl(\Psi(p_1,p_2,s)\bigr).
\]
Since $\mathfrak D_1\subset\mathfrak X(\R^{n_1})$, one has
$(\Phi_2)_*(\mathfrak D_1)=0$, so $\Phi_2$ is independent of $p_1$.
Similarly, $\Phi_1$ is independent of $p_2$. Thus
\[
 \Phi_1=\Phi_1(p_1,s),
 \qquad
 \Phi_2=\Phi_2(p_2,s).
\]

Because $\Theta$ is constant, one has $AT=0$ and
$\overbar\nabla_TT=(\Theta I-P)AT=0$. Hence the integral curves of $T$ are
geodesics in both $\Sigma$ and the ambient product, and each component
$\Phi_i(p_i,s)$ is a geodesic in its factor. Therefore, for
$i\in\{1,2\}$, there is a parallel family
$f_t^i:M_i\to\mathbb Q_{\epsilon_i}^{n_i}$ such that
\[
 \Phi_i(p_i,s)=f_{t_i(s)}^i(p_i).
\]

We may assume, without loss of generality, that
$t_1(s)=s$ and $t_2(s)=\phi(s)$. The constant angle condition makes
the speeds of both component geodesics constant, and proper splitting
makes them nonzero. Thus $\phi(s)=as+b$ with $a\neq0$. Consequently,
we can write
\begin{equation}\label{eq:parametrization}
 \Phi(p_1,p_2,s)
 =\bigl(f_s^1(p_1),f_{as+b}^2(p_2)\bigr),
\end{equation}
that is, locally, $\Sigma$ is an $(f_s^1,f_{\phi(s)}^2)$-hypersurface
with $\phi$ affine. Moreover, since $\Sigma$ has constant mean curvature,
it follows from  Lemma~\ref{lem-Hparallels} that both families
$f_s^1$ and $f_{\phi(s)}^2$ are  isoparametric.

The unit normal vector field of $\Sigma$ is
\begin{equation}\label{eq:normal}
 N=\frac{-a\eta_1+\eta_2}{\mathcal W},
 \qquad
 \mathcal W=\sqrt{1+a^2},
\end{equation}
where $\eta_1$ and $\eta_2$ are the unit normals to $f_s^1$ and
$f_{\phi(s)}^2$, respectively.

For a small $\delta>0$, the map
\[
 \Psi_r(p_1,p_2,s)
 =\exp_{\Phi(p_1,p_2,s)}(rN),
 \qquad r\in J=(-\delta,\delta),
\]
parameterizes the parallel hypersurface to $\Sigma$ at signed distance
$r$. Setting
\[
 \alpha=-\frac a{\mathcal W}\,,
 \qquad
 \beta=\frac1{\mathcal W}\,,
\]
we obtain from~\eqref{eq:parametrization} and~\eqref{eq:normal}
\[
 \Psi_r(p_1,p_2,s)
 =\bigl(f_{s+\alpha r}^1(p_1),
        f_{as+b+\beta r}^2(p_2)\bigr).
\]

Let $H_i(t)$ denote the constant  mean curvature of the
hypersurface $f_t^i$. Using formula~\eqref{eq-H},
the mean curvature of $\Sigma_r$ is
\begin{equation}\label{eq:mean-curvature-parallel}
 H_r(s)
 =\alpha H_1(s+\alpha r)
  +\beta H_2(as+b+\beta r).
\end{equation}

However, because $\Sigma$ is isoparametric, $H_r$ is constant on $\Sigma_r$ for
every sufficiently small $r$.
In particular,
\begin{equation}\label{eq:derivative-H}
 \frac{d}{ds}H_r(s)=0.
\end{equation}

Since $a\beta=-\alpha\neq0$, differentiating
\eqref{eq:mean-curvature-parallel} and using
\eqref{eq:derivative-H} yields
\[
 H_1'(s+\alpha r)=H_2'(as+b+\beta r)
 \qquad\forall (s,r)\in I\times J.
\]
Introducing the variables
\[
 x=s+\alpha r,
 \qquad
 y=as+b+\beta r,
\]
the preceding identity becomes
\begin{equation}\label{eq:H-primes}
 H_1'(x)=H_2'(y).
\end{equation}

Now note that the map $(s,r)\mapsto(x,y)$ is a local diffeomorphism, because its
Jacobian determinant is
\[
 \beta-a\alpha=\mathcal W\neq0.
\]
Hence $x$ and $y$ vary independently in an open set. It follows
then from
\eqref{eq:H-primes} that both sides are constant:
\[
 H_1'(s+\alpha r)=H_2'(as+b+\beta r)=C
 \qquad\forall(s,r)\in I\times J.
\]

It remains to inspect the standard isoparametric families in the two
factors. In Euclidean space they are hyperplanes, spheres, and tubes
about affine subspaces. For a sphere or a non-flat tube, each nonzero
principal curvature of a parallel hypersurface has the form
$$\kappa(t)=\frac{\kappa_0}{1-t\kappa_0},$$
up to the choice of orientation;
hence $H_1'(t)$ is not constant. Therefore $f_t^1$ is a family of
hyperplanes and $H_1'\equiv0$. In hyperbolic space, the non-horospherical
families have principal curvatures involving
$\delta\tanh(\delta t)$ or $\delta\coth(\delta t)$, where $\delta=\sqrt{|\epsilon_2|}$.
Their mean curvature
derivatives are not constant. Thus $f_t^2$ is a family of horospheres
and $H_2'\equiv0$. Therefore $\Sigma$ is an open set of a
flat-horospherical hypersurface.
\end{proof}


Now, we restate and prove Theorem~\ref{thm-classification01}.

\classificationa*
\begin{proof}
\noindent
{\textup{(i)} ${\Rightarrow}$ \textup{(ii)}.}
Suppose that $\Sigma$ is isoparametric. Then, by Theorem~\ref{thm-isoparametric-->thetaconstant},
$\theta$ is constant on $\Sigma$.
Since, in space forms,
any isoparametric hypersurface is an open set of a complete one,
if $\Theta^2\equiv1$,
Proposition~\ref{prop-Ni=0} gives that $\Sigma$ is an open set of a
hypersurface of type \textup{(a)} or \textup{(b)}.
Assume, therefore, that $\Theta^2<1$. If $\epsilon_i>0$ for some
$i\in\{1,2\}$, Proposition~\ref{prop:theta-one}\textup{(i)} would imply
that $\Theta^2\equiv1$, which is a contradiction. Thus
$\epsilon_1,\epsilon_2\leq0$. Since
$\epsilon_1\epsilon_2\leq0$
 and $\epsilon_1^2+\epsilon_2^2\neq0$,
it follows, up to interchanging the factors, that
$\qq=\R^{n_1}\times\h_{\epsilon_2}^{n_2}.$

If $\mathcal A_2$ were nonzero at some point of $\Sigma$, then
Proposition~\ref{prop:theta-one}\textup{(ii)} would again give
$\Theta^2\equiv1$. Hence
$\mathcal A_2\equiv0.$
Since $\Theta^2<1$, we have $T\neq0$, and the equality
$\mathcal A_2\equiv0$ means that the principal directions of $\Sigma$
split properly everywhere. Proposition~\ref{prop:flat-horospherical}
then shows that $\Sigma$ is an open set of a flat-horospherical
hypersurface.

\medskip
\noindent
{\textup{(ii)} ${\Rightarrow}$ \textup{(i)}.}
We have already verified that all the hypersurfaces listed in
\textup{(ii)} are isoparametric and have constant principal
curvatures. This proves \textup{(ii)} $\Rightarrow$ \textup{(i)} and
also shows that every isoparametric hypersurface of $\qq$ has constant
principal curvatures.

\medskip
\noindent
{\textup{(ii)} ${\Rightarrow}$ \textup{(iii)}
\emph{when} ${\epsilon_1,\epsilon_2\leq0}$.}
In this case, every complete isoparametric hypersurface of
$\q_{\epsilon_i}^{n_i}$ is homogeneous. Consequently, the
hypersurfaces of type \textup{(a)} and \textup{(b)} are homogeneous.
We have also verified that flat-horospherical hypersurfaces are
homogeneous. Therefore every hypersurface listed in \textup{(ii)} is
homogeneous in $\qq$, and hence $\Sigma$ satisfies \textup{(iii)}.

\medskip
\noindent
{\textup{(iii)} ${\Rightarrow}$ \textup{(i)}.}
As observed earlier, every homogeneous hypersurface of a Riemannian
manifold is isoparametric. Hence every open subset of a homogeneous
hypersurface is isoparametric.

\medskip
\noindent
{\textup{(iii)} ${\Leftrightarrow}$
\textup{(iv)} \emph{when some} ${\epsilon_i>0}$.}
We have  that \textup{(iv)} implies \textup{(iii)}, because
the product of a homogeneous hypersurface with a complete space form
is homogeneous.

Conversely, suppose that $\Sigma$ is an open set of a homogeneous
hypersurface $$\widehat\Sigma\subset\qq.$$ Then $\widehat\Sigma$ is
isoparametric. Since one of the factors is spherical,
Proposition~\ref{prop:theta-one}\textup{(i)} gives
$\Theta^2\equiv1.$
Proposition~\ref{prop-Ni=0} therefore shows that
$\widehat\Sigma$ is a hypersurface of type \textup{(a)} or
\textup{(b)}. Moreover, the second assertion of that proposition
implies that the corresponding hypersurface
$\Sigma_i^{n_i-1}\subset\q_{\epsilon_i}^{n_i}$ is homogeneous.
Thus \textup{(iii)} implies \textup{(iv)}.
\end{proof}

\begin{remark} \label{rem:gaps02}
In~\cite{tan02}, the authors claim to have classified the
complete isoparametric hypersurfaces of $\s^n\times\R^m$ and $\h^n\times\R^m$.
However, as with their proof of the constancy of the angle, the proof of their
classification theorem (Theorem~1.4) contains several serious gaps. First, it
relies on the constancy of the angle, whose proof is already compromised.
Second, the classification depends heavily on Proposition~3.1, whose proof
contains unjustified arguments and incorrect claims. For instance, in their
rank argument, the authors fail to account for the fact that the focal map
of the full hypersurface and the focal maps of the slices have different
domains and different dimensions. Hence, the identities involving exponential
maps do not imply that their differentials have the same rank. In particular, the map
$$
 p\longmapsto
 \exp_p\left(\frac{2t}{1-C}(0,N^v_p)\right)
$$
is not the normal exponential map of a fixed slice $\Sigma_x$, as claimed,
because its horizontal coordinate also varies with $p$. Consequently,
the asserted synchronization of the focal points of the two slices has not
been established, and this gap affects several essential steps in the classification.
\end{remark}


\section{Proof of Theorem~\ref{thm-classification02}} \label{sec-classification02}

In this final section, we prove Theorem~\ref{thm-classification02}. We first show that, under
the given conditions, the block $\mathcal A_2$ of $A$ has constant norm. Then we establish
an analogue of Proposition~\ref{prop:flat-horospherical}: every connected isoparametric hypersurface of
$\h_{\epsilon_1}^{n_1}\times\h_{\epsilon_2}^{n_2}$
whose principal directions split properly at every point is an open set of a bi-horospherical hypersurface.

\begin{proposition} \label{prop:|A2|}
Assume that $\epsilon_1\epsilon_2\ne 0$ and let $\Sigma\subset\qq$
be an isoparametric hypersurface with constant angle function
$\Theta\in(-1,1)$. Assume further that
$(\sqrt{\tau_1},\sqrt{\tau_2})$
is $(n_1,n_2)$-independent.
Then $\|\mathcal A_2\|$ is constant on $\Sigma$. Moreover, if
$\mathcal A_2=0$ at some point of $\Sigma$, then
$\epsilon_1$ and  $\epsilon_2$ are both negative.
\end{proposition}

\begin{proof}
Since $\theta$ is constant and
$\epsilon_1\epsilon_2\ne 0$,
the function $D(r)$ defined in~\eqref{eq-D&cs} becomes
\begin{equation} \label{eq:drp}
D(r)=\sum_{u=0}^{n_1-1}\sum_{v=0}^{n_2-1}
\alpha_{u,v}(p)\,f_{u,v}(r),
\end{equation}
where
$$
f_{u,v}(r)=c_{\tau_1}^{\,n_1-1-u}(r)s_{\tau_1}^{\,u}(r)
c_{\tau_2}^{\,n_2-1-v}(r)s_{\tau_2}^{\,v}(r).
$$

As we have seen, $D$ is independent of
$p$.  Hence, the equality~\eqref{eq:drp}
applied to $p,q\in\Sigma$ gives
\[
\sum_{u=0}^{n_1-1}\sum_{v=0}^{n_2-1}
\bigl(\alpha_{u,v}(p)-\alpha_{u,v}(q)\bigr)f_{u,v}(r)=0.
\]
But since $\bigl(\sqrt{\tau_1},\sqrt{\tau_2}\bigr)$
is $(n_1,n_2)$-independent, Lemma~\ref{lemma-n1n2dep} implies
that $\{f_{u,v}\}$ are linearly independent. Therefore
\[
\alpha_{u,v}(p)=\alpha_{u,v}(q) \qquad \forall p,q\in\Sigma,
\]
that is, every coefficient  $\alpha_{u,v}$ is constant on $\Sigma$.

Because $\theta$ is constant, we have $AT=0$, which implies that
$$a_{1j}=a_{j1}=0\qquad \forall j\in\{1,\dots, n-1\}.$$
As a consequence,
the coefficient $\alpha_{1,0}$ of $\cone^{n_1-2}(r)\sone(r)\ctwo^{n_2-1}(r)$ is
\begin{equation} \label{eq-alpha100}
\alpha_{1,0}=-\sum_{i=2}^{n_1}a_{ii}=-{\rm tr}\,\mathcal A_1.
\end{equation}

Indeed, by multilinearity of the determinant, $\alpha_{1,0}$
is the sum, over $i\in\{2,\ldots,n_1\}$, of the determinants of the
diagonal matrices $M_i$ obtained
by selecting the diagonal term $-a_{ii}s_{\tau_1}(r)$ in the $i$-th row,
the diagonal terms $c_{\tau_1}(r)$ in the
other rows of the first block, and the diagonal terms $c_{\tau_2}(r)$
in the second block, that is,
$$M_i={\rm diag}\,(\cone(r),\dots,\cone(r),-a_{ii}\sone(r), \cone(r),\dots, \cone(r), \ctwo(r),\dots, \ctwo(r)),$$
with $n_1-2$ occurrences of $\cone(r)$ and $n_2-1$ occurrences of $\ctwo(r)$. Hence
\[
\det M_i=-a_{ii}\cone^{n_1-2}(r)\sone(r)\ctwo^{n_2-1}(r).
\]

Analogously, for $n_1>2$, we have that
\[
\alpha_{2,0}=\sum_{2\le i<j\le n_1}a_{ii}a_{jj}-a_{ij}^2,
\]
and so (consider $\alpha_{2,0}=0$ if $n_1=2$),
\begin{equation} \label{eq-normA1}
\|\mathcal A_1\|^2=\sum_{i=2}^{n_1}a_{ii}^2+2\sum_{i<j=2}^{n_1}a_{ij}^2=\alpha_{1,0}^2-2\sum_{i<j=2}^{n_1}a_{ii}a_{jj}
+2\sum_{i<j=2}^{n_1}a_{ij}^2=\alpha_{1,0}^2-2\alpha_{2,0}.
\end{equation}

It follows from~\eqref{eq-alpha100}-\eqref{eq-normA1}
that the trace and the norm of
$\mathcal A_1$ are both  constant on $\Sigma$.
In particular,
$${\rm tr}\,\mathcal A_3={\rm tr}\,A-{\rm tr}\,\mathcal A_1=H-{\rm tr}\,\mathcal A_1$$ is constant as well.

Therefore, taking traces on both sides of the
first two equalities in~\eqref{eq:system} yields
\begin{equation}  \label{eq-matricialnorms01}
\left\{
\begin{array}{lcl}
\|\mathcal A_1\|^2 -c\|\mathcal A_2\|^2-(n_1-1)\tau_1&=&0,\\[1ex]
\|\mathcal A_2\|^2-c\|\mathcal A_3\|^2+(n_2-1)c\tau_2&=&0.\\[1ex]
\end{array}
\right.
\end{equation}

As $\|\mathcal A_1\|$, $c$, $n_1$ and $\tau_1$ are constants, the first equality in~\eqref{eq-matricialnorms01}
implies that $\|\mathcal A_2\|$ is constant. Moreover, if $\mathcal A_2=0$, the first equality implies $\tau_1\ge0$, and
the second one  implies $\tau_2\ge 0$. Thus, since we are assuming $\epsilon_1\epsilon_2\ne 0$, we conclude that
$$\epsilon_1=-\tau_1/\|N_1\|^2\qquad\text{and}\qquad\epsilon_2=-\tau_2/\|N_2\|^2$$ are
both negative. This completes the proof.
\end{proof}


\begin{proposition}\label{prop:bihorospherical}
Let
$\Sigma\subset\h_{\epsilon_1}^{n_1}\times\h_{\epsilon_2}^{n_2}$
be a connected isoparametric hypersurface whose principal directions
split properly everywhere. Then $\Sigma$ is an open set of a
bi-horospherical hypersurface.
\end{proposition}

\begin{proof}
The proof will be sketched, since it is analogous to the one given for
Proposition~\ref{prop:flat-horospherical}.

The local-product argument
given there yields parallel families
\[
 f_t^i:M_i\longrightarrow\h_{\epsilon_i}^{n_i},
 \qquad i=1,2,
\]
and constants $a\neq0$ and $b$ such that $\Sigma$ is locally
parameterized by
\[
 \Phi(p_1,p_2,s)=\bigl(f_s^1(p_1),f_{as+b}^2(p_2)\bigr).
\]
Setting
\[
 \mathcal W=\sqrt{1+a^2},
 \qquad
 \alpha=-\frac{a}{\mathcal W},
 \qquad
 \beta=\frac1{\mathcal W},
\]
the parallel hypersurface at signed distance $r$ is locally
parameterized by
\[
 \Phi_r(p_1,p_2,s)=\bigl(f_{s+\alpha r}^1(p_1),
       f_{as+b+\beta r}^2(p_2)\bigr).
\]

As in the proof of Proposition~\ref{prop:flat-horospherical}, the
isoparametric condition first implies that both factor families
$f_t^1$ and $f_t^2$ are isoparametric. Denoting by $H_i(t)$  the
mean curvature of $f_t^i$, the mean curvature of the
parallel hypersurface is
\[
 H_r(s)=
 \alpha H_1(s+\alpha r)
 +\beta H_2(as+b+\beta r).
\]
Since $H_r$ is independent of $s$ and $a\beta=-\alpha\neq0$, we obtain
\[
 H_1'(s+\alpha r)
 =
 H_2'(as+b+\beta r).
\]

Writing
$x=s+\alpha r$ and $y=as+b+\beta r,$
we have that $x$ and $y$ vary independently in an open set, and consequently
\begin{equation}  \label{eq:h'}
 H_1'(x)=H_2'(y)=C
\end{equation}
for some constant $C\in\mathbb R$.

As shown in the proof of Proposition~\ref{prop:flat-horospherical},
the only isoparametric hypersurfaces of $\h_{\epsilon_i}^{n_i}$
satisfying~\eqref{eq:h'} are  the horospheres.
Therefore $\Sigma$ is an open set of a
bi-horospherical hypersurface, as we wished to prove.
\end{proof}


We conclude by restating and proving Theorem~\ref{thm-classification02}.

\classificationb*
\begin{proof}
\noindent
{\textup{(i)} ${\Rightarrow}$ \textup{(ii)}.}
Suppose that $\Sigma$ is isoparametric. By
Theorem~\ref{thm-isoparametric-->thetaconstant}, its angle function
$\Theta$ is constant.
If $\Theta(p)^2=1$, then $\Theta^2\equiv1$ on $\Sigma$, and
Proposition~\ref{prop-Ni=0} shows that $\Sigma$ is an open set of a
hypersurface of type \textup{(a)} or \textup{(b)}.

Assume, therefore, that $\Theta(p)^2<1$. Since $\Theta$ is constant,
we have $\Theta^2<1$
everywhere on $\Sigma$. Moreover, because $p$ is a distinguished
point, the principal directions of $\Sigma$ split properly at $p$;
equivalently, $\mathcal A_2(p)=0.$
In addition, $\tau_1$ and $\tau_2$ are constant. Hence the
independence assumption at $p$ holds everywhere,
and Proposition~\ref{prop:|A2|} applies. Consequently,
$\|\mathcal A_2\|$ is constant on $\Sigma$. Since
$\mathcal A_2(p)=0$, we conclude that
$\mathcal A_2\equiv0.$

The final assertion of Proposition~\ref{prop:|A2|} also gives
$\epsilon_i<0$, and so $\mathbb Q_{\epsilon_i}^{n_i}=\mathbb H_{\epsilon_i}^{n_i}.$
In addition,
the equality $\mathcal A_2\equiv0$ means that the principal directions
of $\Sigma$ split properly everywhere. Thus, by
Proposition~\ref{prop:bihorospherical}, $\Sigma$ is an open
set of a bi-horospherical hypersurface. This proves
\textup{(i)} $\Rightarrow$ \textup{(ii)}.

\medskip
\noindent
{\textup{(ii)} ${\Rightarrow}$ \textup{(i)}.}
We have already verified that all the hypersurfaces listed in
\textup{(ii)} are isoparametric and have constant principal
curvatures. This proves \textup{(ii)} $\Rightarrow$ \textup{(i)} and
also shows that every isoparametric hypersurface satisfying the
hypotheses of the theorem has constant principal curvatures.

\medskip
\noindent
{\textup{(ii)} ${\Rightarrow}$ \textup{(iii)}
\emph{when} ${\epsilon_1,\epsilon_2<0}$.}
Every complete isoparametric hypersurface of
$\h_{\epsilon_i}^{n_i}$ is homogeneous. Consequently, the
hypersurfaces of type \textup{(a)} and \textup{(b)} are homogeneous.
We have also verified that bi-horospherical hypersurfaces are
homogeneous. Therefore every hypersurface listed in \textup{(ii)} is
homogeneous in $\h_{\epsilon_1}^{n_1}\times\h_{\epsilon_2}^{n_2},$
and hence $\Sigma$ satisfies \textup{(iii)}.

\medskip
\noindent
{\textup{(iii)} ${\Rightarrow}$ \textup{(i)}.}
As observed earlier, every homogeneous hypersurface of a Riemannian
manifold is isoparametric. Hence every open subset of a homogeneous
hypersurface is isoparametric.

\medskip
\noindent
{\textup{(iii)} ${\Leftrightarrow}$
\textup{(iv)} \emph{when} ${\epsilon_1,\epsilon_2>0}$.}
As before, we have  that \textup{(iv)} implies \textup{(iii)} because
the product of a homogeneous hypersurface with a complete space form
is homogeneous.

Conversely, suppose that $\Sigma$ is an open subset of a homogeneous
hypersurface
$\widehat\Sigma\subset\s_{\epsilon_1}^{n_1}\times\s_{\epsilon_2}^{n_2}.$
Then $\widehat\Sigma$ is isoparametric. By the implication
\textup{(i)} $\Rightarrow$ \textup{(ii)}, already proved,
$\widehat\Sigma$ is of type \textup{(a)} or \textup{(b)}, since
bi-horospherical hypersurfaces occur only in products of hyperbolic
spaces. The second assertion of Proposition~\ref{prop-Ni=0} then
implies that the corresponding hypersurface
$\Sigma_i^{n_i-1}\subset\s_{\epsilon_i}^{n_i}$
is homogeneous. Thus \textup{(iii)} implies \textup{(iv)}.
\end{proof}


\bigskip
\noindent
{\bf Acknowledgments.} We are grateful to Jinxuan Chen
for pointing out some inconsistencies in the proof of
Theorem~\ref{thm-isoparametric-->thetaconstant} in a preliminary version
of this paper. Giuseppe Pipoli was partially supported by the INdAM--GNSAGA Project,
CUP E53C24001950001, and by PRIN Project 20225J97H5.


\bigskip
\noindent
{\bf Data availability.}  Not applicable.


\bigskip
\noindent
{\bf Conflict of interest.}  The authors declare that there is no conflict of interest.


\end{document}